\documentclass{stavanger}

\usepackage[utf8]{inputenc}
\usepackage{amssymb,amsmath,amsfonts,amsthm}
\usepackage{ dsfont }
\usepackage{hyperref}
\usepackage{natbib}
\usepackage{graphicx}
\usepackage{tikz}
\usepackage{cleveref}
\usepackage{mathtools}
\usepackage{placeins}
\usepackage[export]{adjustbox} 

\begin{document}


\begin{frontmatter}

\begin{fmbox}


\dochead{Research Article - Preprint}{FP}


\title{A new variational discretization technique for initial value problems bypassing governing equations}


\author[
   addressref={aff1},                   	  
   corref={aff1},                     		  
   email={alexander.rothkopf@uis.no}   		  
]{\inits{AR}\fnm{Alexander} \snm{Rothkopf}}
\author[
   addressref={aff2,aff3},                   	  
   email={jan.nordstrom@liu.se}   		  
]{\inits{JN}\fnm{Jan} \snm{Nordstr{\"o}m}}


\address[id=aff1]{
  \orgname{Faculty of Science and Technology}, 	 
  \street{University of Stavanger},                     		 
  \postcode{4021},                               			 
  \city{Stavanger},                              				 
  \cny{Norway}                                   				 
}

\address[id=aff2]{
  \orgname{Department of Mathematics, Applied Mathematics}, 	 
  \street{Link{\"o}ping University},                     		 
  \postcode{SE-581 83},                               			 
  \city{Link{\"o}ping},                              				 
  \cny{Sweden}                                   				 
}

\address[id=aff3]{
  \orgname{Department of Mathematics and Applied Mathematics}, 	 
  \street{University of Johannesburg},                     		 
  \postcode{P.O. Box 524, Auckland Park 2006},                           
  \city{Johannesburg},                              				 
  \cny{South Africa}                                   				 
}



\end{fmbox}


\begin{abstractbox}

\begin{abstract} 
Motivated by the fact that both the classical and quantum description of nature rest on causality and a variational principle, we develop a novel and highly versatile discretization prescription for classical initial value problems (IVPs). It is based on an optimization (action) functional with doubled degrees of freedom, which is discretized using a single regularized summation-by-parts (SBP) operator. Formulated as optimization task it allows us to obtain classical trajectories without the need to derive an equation of motion.
The novel regularization we develop in this context is inspired by the weak imposition of initial data, often deployed in the modern treatment of IVPs and is implemented using affine coordinates. We demonstrate numerically the stability, accuracy and convergence properties of our approach in systems with classical equations of motion featuring both first and second order derivatives in time.
onvergence properties of our approach in systems with classical equations of motion featuring both first and second order derivatives in time.
\end{abstract}



\end{abstractbox}


\end{frontmatter}


\section{Introduction}

The numerical treatment of dynamical phenomena in classical and quantum systems is at the core of progress in natural sciences and engineering. In computational fluid dynamics \cite{blazek2015computational} or electrodynamics \cite{taflove2005computational}, a set of coupled partial differential equations is solved on a predefined geometric domain with boundary conditions, starting from an initial condition, in order to predict trajectories of point particles or configurations of fields. In the study of atomic properties, linear and non-linear variants of the Schr\"odinger equation or more generally Lindblad equations \cite{breuer2002theory} of multiple entangled particles are solved as initial value problems. For an understanding of the nuclei of atoms on the other hand, an ensemble of fluctuating quantum fields of a non-linear variant of Maxwell's equations (Yang-Mills theory) needs to be simulated on a hypercubic grid (lattice QCD) \cite{gattringer2009quantum}.

Much progress has been made in developing accurate and cost effective discretization schemes for partial differential equations over the past two decades. Due to their ease of implementation, finite difference schemes have long enjoyed popularity, but historically were challenged when confronted with intricate simulation geometries. It took the development of summation-by-parts (SBP) finite difference operators (for reviews see e.g. \cite{svard2014review,fernandez2014review,lundquist2014sbp}), to elevate finite difference schemes to a similar level of versatility as traditional functional basis approaches, such as Galerkin schemes \cite{nordstrom2017roadmap}. The SBP approach both in spatial dimensions, as well as in time \cite{lundquist2014sbp,nordstrom2013summation,nordstrom2016summation} provides proofs of stability for finite difference based discretization schemes via the so-called energy method and is easily extended to higher order approximations. 

Implementing the integration-by-parts property of the underlying continuum IVP, summation-by-parts operators are an example of so-called mimetic discretizations. It has been shown that SBP operators form a versatile framework, which encompasses various other numerical approximation techniques besides finite differences \cite{svard_stable_2007}, such as finite volume schemes \cite{nordstrom_finite_2003}, spectral element \cite{carpenter_entropy_2014}, flux reconstruction \cite{ranocha_summation-by-parts_2016} and both continuous \cite{abgrall_analysis_2020} and discontinuous Galerkin (dG) \cite{gassner_skew-symmetric_2013,hesthaven_stable_1996} schemes.

A crucial development in the numerical treatment of differential equations is the concept of weak boundary or initial conditions. It acknowledges that the solution of a discretized PDE  not only in the interior of the domain but also on the boundary ( or initial time slice ) need only be as accurate as the order of the discretization. By allowing the solution to deviate from the initial or boundary conditions within the tolerance of the discretization, one obtains a new lever, which one can exploit in the construction of discretization schemes. The simultaneous approximation term (SAT) approach \cite{carpenter1994time} e.g. implements weak boundary or initial conditions by the addition of appropriately designed penalty terms to the differential equation of interest. In recent studies it has been shown how to absorb part of these penalty terms into a redefinition of the SBP operators, in order to reduce their null-space to the corresponding physical dimension, leading to so called null-space consistent SBP operators \cite{svard2019convergence,linders_properties_2020,svard_convergence_2021,ranocha2021new}.

In spite of this substantial progress in the numerical treatment of IVPs, challenges of both conceptual and technical nature remain. The treatment of intrinsic constraints, e.g. the divergence constraint for Maxwell's equations and the discretization of second order systems are two examples. The third one is the derivation of the equations of motion of linear electrodynamics or non-linear Yang-Mills theory in terms of the gauge potentials (see e.g. \cite{Ipp:2018hai}). It is relevant for quantum theory and depends on a choice of gauge. The choice of setting the zeroth component of the four-potential to zero $A_0=0$ renders the role of Gauss' law opaque, since $A_0$ is actually the Lagrange multiplier that preserves this intrinsic constraint. On the other hand, the discretized Lagrangian of these theories remains manifestly gauge invariant and $A_0$ retains its central role. If one could solve the associated initial value problem, i.e. determine future field configurations directly on the level of the action, without the need to derive the equation of motion, no choice of gauge is necessary and the manifested gauge invariance would render Gauss' law automatically fulfilled.

In the treatment of initial value problems for second order ODEs with the SBP-SAT technique it was found that using the same regularized SBP operator for the first and second derivative does not lead to stable procedures. Instead, different SBP operators for the position and velocity degree of freedom had to be defined \cite{nordstrom2016summation}. On the other hand in the action formulation of second order systems, at most first order derivatives act on the degrees of freedom \cite{arnold_mathematical_1989}. This reduction of the order of the derivatives compared to the corresponding governing equations by at least one power is a general feature of the action formulation. A lower order derivative operator is more robust against rough data than a higher order one, which requires more regularity \cite{gustafsson2007high}. By expressing the theory directly on the level of the action one can reduce the highest order of derivatives and in particular for second order systems, a single modified SBP operator suffices to obtain a regularized and null-space consistent discretization.

In this paper we develop a discretization scheme for one-dimensional initial value problems, based on a generalized variational principle applied directly to the action of the system of interest. I.e. we will determine the classical trajectory of the system, given an initial condition, without the need to derive an equation of motion. To this end we deploy the SBP technique to approximate derivatives occurring in the continuum formulation of the problem and will take inspiration by the SAT approach to regularize the resulting difference operators. This is achieved by the use of affine coordinates to absorb the whole penalty term including data into a redefined null-space consistent SBP operator.

Of course all examples treated in this study have well known ODEs as equation of motion, which can be solved with established numerical methods. As proof of principle, the intention of this study is merely to establish the validity of the direct action based approach and prepare the ground for its application to higher dimensional theories in future work, where its utility is more obvious. Our direct action approach allows us to formulate the system using any geometry for which SBP operators have been developed and therefore allows us to avoid the intricate construction of dual meshes \cite{elcott_building_2005}, which underlie implementations of other variational approaches to IVPs, such as discrete exterior calculus \cite{stern_geometric_2015}.

It is important to note that the main point of this paper is to present a new technique for solving IVPs \textit{without utilizing governing equations}.

The paper is structured as follows: In \cref{sec:contForm1} we review the continuum formulation of the conventional variational principle of classical physics for initial value problems, featuring equations of motion with second order derivatives in time. In the subsequent \cref{sec:SBPdiscr} we introduce our discretization prescription, discuss the need for regularization and construct a regularized SBP operator based on initial value data. To extend the applicability of our discretization scheme to a more general class of systems, we consider a generalized variational principle in \cref{sec:generalizedVP} and show that our approach successfully captures differential equations of motion containing also first order derivatives. We close with a brief summary in \cref{sec:summary}.

\section{Continuum formalism for second order equations of motion}
\label{sec:contForm1}

The classical physics of closed systems (i.e. systems that are not in contact with their environment) is conveniently captured via their Lagrangian. In point mechanics the Lagrangian is a functional, which depends on the trajectory of the point mass $x(t)$ and its velocity $\dot x(t)$. In a field theory, such as in electromagnetism, it is formulated in terms of the vector fields $A_\mu(x)$ and their derivatives $\partial_\nu A_\mu(x)$. In the simple systems under consideration here, the Lagrangian can be written as the difference between the kinetic energy of the system and its potential energy. Taking a point mass in a constant gravitational field as an explicit example we have
\begin{align}
    {\cal L}=T-V= \frac{1}{2}m\dot x^2(t) - mgx(t).
\end{align}

In the 19th century, physicists discovered that the trajectory of a classical particle can be obtained from a variational principle (see e.g. ref.\cite{goldstein_classical_2001}). If a point mass starts out at position $x(t_1)$ at $t_1$ and ends up at position $x(t_2)$ at $t_2$ then the classical trajectory that is realized in nature between those points is given by the critical point of the classical action\footnote{More specifically, for short times, where the classical trajectory has not yet reached any turning point, the action exhibits an actual minimum. In the presence of turning points in the classical trajectory it in general represents a saddle point of the action. (see e.g. ref.~\cite{arnold_mathematical_1989})}
\begin{align}
    S[x(t),\dot x(t)]=\int_{t_1}^{t_2} dt \, {\cal L}[x(t),\dot x(t)] .
\end{align}
This observation is stated as a \textit{boundary value problem}, where the start and end point of the trajectory are specified. While it is of conceptual interest that such a formulation exists, we cannot use it to determine the classical trajectory itself in a causal fashion, since in order to formulate the variational principle, we already need to know where the point mass will end up at $t_2$. When setting up an experiment, we are of course only in control of the initial position and velocity.

To proceed, physicists conventionally convert the above boundary value problem into an \textit{initial value problem} using the following strategy: one derives a set of differential equations that are equivalent to the variational principle and which can be solved as initial value problem. These are the celebrated Euler-Lagrange equations. We wish to inspect the variation of the Lagrangian using a slightly deformed path $x(t)+\delta x(t)$. Here $\delta x(t)$ is an arbitrary function with the only condition that it goes to zero at $t_1$ and $t_2$, as the points $x(t_1)$ and $x(t_2)$ are fixed. Varying the action we obtain
\begin{align}
    \delta {\cal S} &= \int_{t_1}^{t_2} dt \Big\{ \frac{\partial {\cal L}}{\partial x} \delta x + \frac{\partial {\cal L}}{\partial \dot x} \delta \dot x \Big\}
    =\int_{t_1}^{t_2} dt \Big\{ \frac{\partial {\cal L}}{\partial x} \delta x + \frac{\partial {\cal L}}{\partial \dot x} \frac{d}{dt} \delta x \Big\}\label{eq:equivELstep2}   \\
     &=\int_{t_1}^{t_2} dt \Big\{ \frac{\partial {\cal L}}{\partial x}- \frac{d}{dt}\frac{\partial {\cal L}}{\partial \dot x}   \Big\}\delta x + \left.\Big[ \frac{\partial {\cal L}}{\partial \dot x} \delta x \Big]\right|_{t_1}^{t_2}\label{eq:equivELBVP},
\end{align}
where the last line results from \textit{integration by parts} (IBP). Since the variation $\delta x$ by construction vanishes on the boundary, the term in the square brackets also vanishes. If we inspect the critical point of the functional ${\cal S}$, defined by $\delta {\cal S} =0$, we find that it is equivalent to the term in the curly brackets equalling zero, since $\delta x$ can be any (well behaved) function between $t_1$ and $t_2$.

\begin{figure}
    \centering
    \includegraphics[scale=0.5]{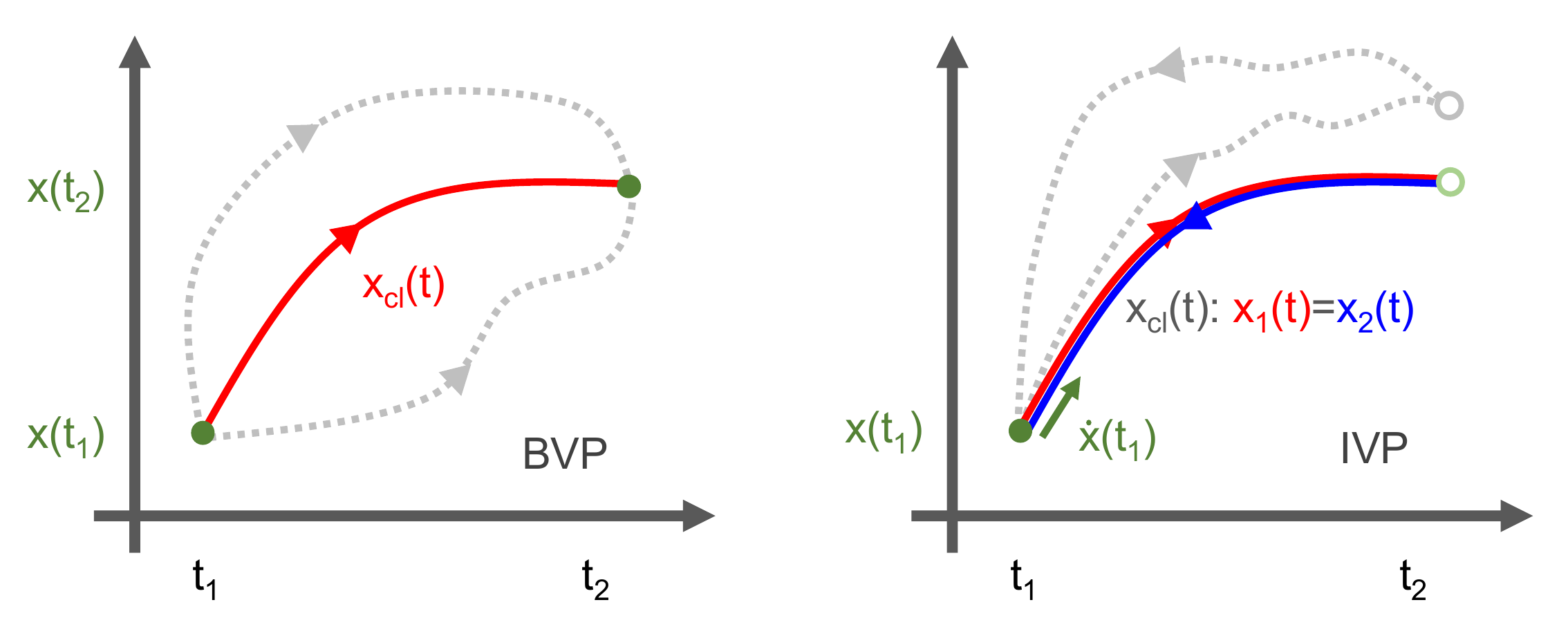}
    \caption{Differences between the variational principle as boundary value problem (left) and as initial value problem (right). In the BVP setting there may exist multiple paths that fulfill the boundary conditions, but only one of them, $x_{\rm cl}(t)$, represents an extremum of the action. In the IVP setting, a doubling of the degrees of freedom is required in order to allow the path at time $t_2$ to vary freely. Figure adapted from ref.~\cite{galley_classical_2013}.}
    \label{fig:sketch}
\end{figure}

 In other words, if we assume the validity of the variational principle, i.e. that the classical trajectory follows from the critical point of the action, then this trajectory must fulfill the Euler-Lagrange equations, which are just the terms inside the curly brackets set to zero
 \begin{align}
\left.\frac{\delta {\cal S}[x,\dot x]}{\delta x}\right|_{x=x_{\rm classical}}=0\quad \underset{\rm IVP}{\overset{\rm BVP}{\iff}} \quad \frac{\partial {\cal L}}{\partial x}- \frac{d}{dt}\frac{\partial {\cal L}}{\partial \dot x} =0.\label{eq:varprincBVP}
\end{align} 
Since for more complex systems (with internal constraints etc.) it is often easier to formulate the action than to derive the Euler-Lagrange equations as initial value problem, our goal here is to formulate and solve the initial value problem as a variational problem directly on the level of the action. To this end we follow the reasoning of ref.~\cite{galley_classical_2013}, which establishes the continuum formalism for the variational IVP.

Retracing the train of thought of ref.\cite{galley_classical_2013}, we first note that the equivalence between the Euler-Lagrange equations (which provide the correct classical equations of motion) and the stationarity of the action requires that the variation of the path vanishes at initial $t_1$ and final time $t_2$ (see \cref{eq:equivELBVP}). Since we do not know $x(t_2)$ apriori and we only know $x(t_1)$ and $\dot x(t_1)$, we must instead find a way how to formulate the variational principle in such a way that the value of $x$ at $t_2$ does not need to be fixed. This can be accomplished by doubling the degrees of freedom with one trajectory $x_1(t)$ describing a \textit{forward} path and one trajectory $x_2(t)$ describing a \textit{backward} path. Intuitively we will use the doubled degrees of freedom in a way reminiscent of the shooting method, extended such that the aim is to find the trajectory which returns to the starting point given by the initial conditions (see \cref{fig:sketch} for a sketch of the difference in the approaches).

In order to achieve the necessary cancellation of the boundary terms between the forward and backward path at $t_2$, ref.\cite{galley_classical_2013} constructs a new joint action for the two degrees of freedom $x_1(t)$ and $x_2(t)$ as
\begin{align}
    S_{\rm IVP}[x_1(t),\dot x_1(t),x_2(t),\dot x_2(t)]&=\int_{t_1}^{t_2} dt\Big( {\cal L}[x_1(t),\dot x_1(t)] - {\cal L}[x_2(t),\dot x_2(t)] \Big),\\&=\int_{t_1}^{t_2} dt \Big( {\rm L}[x_1(t),\dot x_1(t),x_2(t),\dot x_2(t)] \Big).
\end{align}
The Lagrangian housing the backward path $x_2$ is introduced with a relative minus sign, which, as we will show, allows the boundary terms arising in the variation of $ {\cal L}[x_1(t),\dot x_1(t)]$ and ${\cal L}[x_2(t),\dot x_2(t)]$ to cancel. Let's carry out the variation of this new $S_{\rm IVP}$ explicitly, which yields twice as many terms
\begin{align}
    \delta {\cal S} &= \int dt \Big( \Big\{ \frac{\partial {\rm L}}{\partial x_1}- \frac{d}{dt}\frac{\partial {\rm L}}{\partial \dot x_1}   \Big\}\delta x_1 -\Big\{ \frac{\partial {\rm L}}{\partial x_2}- \frac{d}{dt}\frac{\partial {\rm L}}{\partial \dot x_2}   \Big\}\delta x_2\Big)\\
    &+ \left. \Big[ \frac{\partial {\rm L}}{\partial \dot x_1} \delta x_1 \Big]\right|_{t_1}^{t_2} - \left.\Big[ \frac{\partial {\rm L}}{\partial \dot x_2} \delta x_2 \Big]\right|_{t_1}^{t_2}  \label{eq:equivEL}.
\end{align}

In order to see how the cancellations come about, it is advantageous to change coordinates, going over to relative $x_-=x_1-x_2$ and centered coordinates $x_+=(x_1+x_2)/2$. This change is not necessary, but expressed in $x_-$ and $x_+$ the new variational principle can be formulated in a very concise form and the relation between the functional and the resulting differential equation of motion becomes much more lucid.

We vary the action using $x_{\pm}(t)+\delta x_{\pm}(t)$. The new path deformations $\delta x_{\pm}(t)$ vanish at the initial time $t_1$, as the original deformations are set to zero there $\delta x_1(t_1)=\delta x_2(t_1)=0$.
As the action is now a functional of the newly introduced paths ${\cal S}_{\rm IVP}[x_+(t),\dot x_+(t),x_-(t),\dot x_-(t)]$, its variation produces the following expression
\begin{align}
   \nonumber  \delta {\cal S}_{\rm IVP} &= \int dt \Big( \Big\{ \frac{\partial {\rm L}}{\partial x_+}- \frac{d}{dt}\frac{\partial {\rm L}}{\partial \dot x_+}   \Big\}\delta x_+ +\Big\{ \frac{\partial {\rm L}}{\partial x_-}- \frac{d}{dt}\frac{\partial {\rm L}}{\partial \dot x_-}   \Big\}\delta x_- \Big)\\ &+
    \left. \Big[  \frac{\delta {\rm L}}{\delta \dot x_+} \delta x_+ + \frac{\delta {\rm L}}{\delta \dot x_-} \delta x_-    \Big]\right|_{t_1}^{t_2}\label{eq:equivELS2}.
\end{align}

In order to correctly cancel the boundary contribution $ \frac{\delta {\rm L}}{\delta \dot x_-} \delta x_-$ at $t_2$, we see that the values of $x_1(t_2)$ and $x_2(t_2)$ have to agree, i.e. $x_-(t_2)=0$. It is important to note that the paths $x_1$ and $x_2$ themselves are not fixed to a certain value at $t_2$, since we do not know that value apriori. I.e. the forward and backward paths need to be connected, corresponding to the condition $x_1(t_2)=x_2(t_2)$.

What happens to the other boundary term $ \frac{\delta {\rm L}}{\delta \dot x_+} \delta x_+$? Since $\dot x_1(t)=\dot x_+(t)+\frac{1}{2}\dot x_-(t)$ and $\dot x_2(t)=\dot x_+(t)-\frac{1}{2}\dot x_-(t)$, we find the following expression for the derivative of the joint Lagrangian ${\rm L}$
\begin{align}
\frac{\delta {\rm L}}{\delta \dot x_+} = \frac{\delta {\rm L}}{\delta \dot x_1}\frac{\partial \dot x_1}{\partial \dot x_+} + \frac{\delta {\rm L}}{\delta \dot x_2}\frac{\partial \dot x_2}{\partial \dot x_+}= \frac{\delta {\cal L}}{\delta \dot x_1} - \frac{\delta {\cal L}}{\delta \dot x_2}=\pi_1-\pi_2.
\end{align}
In the second equality we have explicitly written ${\rm L}$ as the difference between the individual Lagrangians ${\cal L}$ for the forward path $x_1$ and backward path $x_2$. In the last step we furthermore introduced the conjugate momenta of the paths, defined as $\pi_{1,2}=\delta {\cal L}/\delta \dot x_{1,2}$. This relation between the functional derivative with respect to $x_+$ and the difference between the momenta on the forward and backward path tells us that we can make the remaining boundary term in \cref{eq:equivELS2} vanish at $t_2$ if we construct our paths such that the difference between the momenta $\pi_1(t_2)-\pi_2(t_2)=0$ vanishes at time $t_2$.

For the systems considered here, which exhibit second order derivatives in their equation of motion, we have the kinetic term $T=\frac{1}{2} \dot x^2$ in the Lagrangian ${\cal L}$. This term leads to the identification $\pi_{1,2}=\dot x_{1,2}$. In turn we find that if we require that in addition to the values of the paths at $t_2$ also the derivatives are identified $\dot x_1(t_2)=\dot x_2(t_2)$, both boundary terms in \cref{eq:equivELS2} will vanish. This establishes the necessary conditions for joining the forward and backward path
\begin{align}
x_1(t_2)=x_2(t_2), \qquad \dot x_1(t_2)=\dot x_2(t_2),
\end{align}
in order to relate the extremum of the joint functional ${\rm L}$ to the Euler-Lagrange equation expressions in the curly brackets in \cref{eq:equivELS2}.

We had to introduce doubled degrees of freedom to correctly cancel the boundary terms that arise from the fact that for an IVP the value of the classical path is unknown at time $t_2$. In the end there only exists a single classical trajectory and we hence must undo the proliferation of degrees of freedom. To this end ref.~\cite{galley_classical_2013} introduces what they call the \textit{physical limit}, which enforces $x_1(t)-x_2(t)=x_-(t)=0$ at all times. When applied to the equations of motion resulting from  \cref{eq:equivELS2}, i.e.
\begin{align}
\frac{\partial L}{\partial x_\pm}=\frac{d}{dt}\frac{\partial L}{\partial \dot x_\pm}\label{eq:ELdbl}
\end{align}
only those equations independent of $x_-$ survive. Since the functional $L={\cal L}[x_1,\dot x_1]- {\cal L}[x_2,\dot x_2]$ is constructed from a difference of the Lagrangians on the forward and backward paths it will always contain at least a linear dependence on $x_-$ and $\dot x_-$. Thus only the equation in \cref{eq:ELdbl}, in which the derivative with respect to $x_-$ is taken can survive.
 
Combining the variation of the joint action of the forward and backward path with the \textit{physical limit}, we thus arrive at the following concise formulation of the variational principle for a classical initial value problem 
\begin{align}
\left. \frac{\delta S_{\rm IVP}[x_{\pm}]}{\delta x_-}\right|_{x_-=0,x_+=x_{\rm class}}=0.\label{eq:varIVPdefEq}
\end{align}
Note that in deriving \cref{eq:varIVPdefEq}, integration by parts (IBP) took center stage. This fact motivates the use of summation-by-parts (SBP) operators in the discretization of the variational principle in the next \cref{sec:SBPdiscr}. Indeed, if the discretization is able to exactly mimic IBP, all steps up to this point follow through also in the discrete setting (see \cref{eq:discrEOMderiv} in \cref{sec:firstordersys}).

Formulating classical mechanics as variational problem offers further insight derived from Noether's theorem. Following ref.~\cite{sieberer2016keldysh} one can show that Noether's theorem for an action with doubled degrees of freedom can be established and it provides two important results. Using as starting point the action ${\cal S}_{\rm IVP}$ and using only integration by parts and the swapping of differentiation and variation, it follows that the sum of the energy of the forward and the backward path is preserved in time, as is the difference between the two. This establishes that even though the additional backward path has been added to the system the energy associated with it remains bounded and the system is in fact stable. For more details and discussion see \ref{sec:Noether}.

Take as explicit example the point mass in a constant gravitational field. Its Lagrangian is ${\cal L}=\frac{1}{2}m\dot x^2(t) - mgx(t)$ and the Euler-Lagrange equation reads
\begin{align}
    \ddot x_{\rm class}(t)=-g, \qquad x_{\rm class}(t)=-\frac{1}{2}gt^2 + \dot x(0)t + x(0), \label{eq:ELpointmass}
\end{align}
which is nothing but Newtons law in terms of acceleration and can be solved in a straight forward manner. We will take $g$ to be positive to indicate that gravity is acting downwards. 

Using the formalism based on the doubled degrees of freedom we have instead
\begin{align}
    S_{\rm IVP}&=\int dt\, \Big( \frac{1}{2}m (\dot x^2_1(t)-\dot x^2_2(t)) - mg( x_1(t)-x_2(t))\Big), \\
    &= \int dt \Big( m \dot x_+(t)\dot x_-(t) - mgx_-(t)\Big).\label{eq:exampleconstr}
\end{align}
In computing the variation of the action, we carried out one integration by parts, which, in effect, allows us to re-express $\delta {\cal S}$ as depending solely on the variation of the paths and not on their derivatives (see \cref{eq:equivELstep2}). Similarly we can integrate by parts here to move the time derivative on $x_-$ in the kinetic term to $x_+$
\begin{align}
    S_{\rm IVP}= \int dt \Big(- m \ddot x_+(t) x_-(t) - mgx_-(t)\Big).\label{eq:exampleconstr2}
\end{align}
Since we identify both the values and derivatives of the paths at $t_2$ no boundary terms contribute. Taking the functional derivative of \cref{eq:exampleconstr2} with respect to $x_-$, setting the result to zero and identifying $x_+=x_{\rm class}$ in the \textit{physical limit} yields exactly the conventional Euler-Lagrange equation
\begin{align}
    \ddot x_+(t)=\ddot x_{\rm class}(t)=-g.\label{eq:classeom}
\end{align}

We have by now seen how the continuum variational principle for IVPs is derived and have acquired intuition in a simple system what form the joint action ${\cal S}_{\rm IVP}$ takes on in terms of $x_+$ and $x_-$. In the remainder of the paper we will only work on the level of the joint action and not need to refer to the equation of motion anymore. Let us briefly mention that the corresponding functional $L[x_+,\dot x_+,x_-,\dot x_-]$ for a large variety of systems with second order equations of motion of the form $\ddot x + f(x) =0$ can be written as 
\begin{align}
{\cal S}_{\rm IVP}=\int dt \Big( m \dot x_+(t)\dot x_-(t) - f(x_+)x_-(t)\Big)\label{eq:constrsecord}.
\end{align}

In \cref{sec:generalizedVP}, after having established the discrete formalism for systems with a second order equation of motion in time, we will consider a generalized variational principle also derived in ref. \cite{galley_classical_2013}, which will allow us to extend the discrete treatment to systems with differential equations of motion containing also single derivatives in time.

\section{Variational IVP based on SBP operators in time}
\label{sec:SBPdiscr}
\vspace{0.2cm}
\subsection{A naive SBP discretization of the model boundary value problem}
\label{sec:discbvp}
As a first step, let us formulate the discretized variational principle in its conventional form as boundary value problem. The point mass in a constant gravitational field will again serve as an explicit example. 

Our goal here is to discretize the action for the single forward path $x(t)$
\begin{align}
{\cal S}=\int dt\, \Big( \frac{1}{2}m\dot x^2(t)-mgx(t) \Big),
\end{align}
with Dirichlet boundary conditions $x(0)=x_i$ and $x(1)=x_f$, in order to compute the classical trajectory at its critical point. To this end we introduce the path ${\bf x}=(x(0),x(\Delta t),x(2\Delta t),\ldots)^{\rm T}$ resolved at $N_t$ points with time step $\Delta t =1/(N_t-1)$. The integral can be approximated by a quadrature rule, whose particular form is captured in a (diagonal) positive definite matrix $\mathds{H}$ and which defines an inner product on discretized paths $({\bf x},{\bf x}^\prime)={\bf x}^{\rm T} \mathds{H} {\bf x}^\prime$. 

Remember that the derivation of the continuum variational principle in \cref{eq:equivELBVP} and \cref{eq:equivEL} required us to carry out \textit{integration by parts}. In order to guarantee the equivalence between the stationarity of the action and the solution of the Euler-Lagrange equation, we must therefore use a discretization that mimics IBP in the discrete setting. Hence we discretize the time derivative with a summation by parts operator $\mathds{D}=\mathds{H}^{-1}\mathds{Q}$, compatible with our choice of $\mathds{H}$, where $\mathds{Q}^{\rm T}+\mathds{Q}=\mathds{E}_N-\mathds{E}_0={\rm diag}[-1,0,\ldots,0,1]$.

The lowest order SBP discretization scheme \texttt{SBP21} of order two in the interior and order one on the boundary ensues when choosing the trapezoid rule for integration
\begin{equation}
\nonumber \mathds{H}^{[2,1]}=\Delta t \left[ \begin{array}{ccccc} 1/2 & & & & \\ &1 & & &\\ & &\ddots && \\ &&&1&\\ &&&&1/2 \end{array} \right],
\quad 
\mathds{D}^{[2,1]}=
\frac{1}{2 \Delta t}
\left[ \begin{array}{ccccc} -2 &2 & & &\\ -1& 0& 1& &\\ & &\ddots && \\ &&-1&0&1\\ &&&-2&2 \end{array} \right].
\end{equation}
The next higher order SBP scheme \texttt{SBP42} is fourth order in the interior and second order on the boundary
\begin{align}
\nonumber &\mathds{H}^{[4,2]}=\Delta t \left[ \begin{array}{ccccccc} 
\frac{17}{48} & & & & & & \\
 & \frac{59}{48} & & & & & \\
 & &\frac{43}{48} & & & & \\
 & & &\frac{49}{48} & & & \\
 & & & & &1 & \\
 & & & & & &\ddots \\
 \end{array} \right]\end{align}

 \begin{align}
&\mathds{D}^{[4,2]}=
\frac{1}{\Delta t}
\left[ \begin{array}{ccccccccc} 
-\frac{24}{17}&\frac{59}{34} & -\frac{4}{17} & -\frac{3}{34} & && & & \\
-\frac{1}{2}& 0 & \frac{1}{2} & 0 & && & & \\
 \frac{4}{43}& -\frac{59}{86} & 0 & \frac{59}{86}&-\frac{4}{43} && & & \\
\frac{3}{98}& 0& -\frac{59}{86}  & 0&\frac{32}{49}&-\frac{4}{49}& & & \\
& &\frac{1}{12}  & -\frac{2}{3}&0&\frac{2}{3}& -\frac{1}{12}& & \\
&&&&&&&\ddots&
 \end{array} \right].
\end{align}
In this section we will show explicit results based on the \texttt{SBP21} operator and include the outcomes from the \texttt{SBP42} operator in our scaling tests.

The discretized action, on which the variational boundary value problem rests, reads
\begin{align}
\mathds{S}_{\rm BVP}= \frac{1}{2}m \big(\mathds{D}{\bf x}\big)^{\rm T} \mathds{H} \big(\mathds{D}{\bf x}\big) - mg \mathds{1}^{\rm T} \mathds{H} {\bf x} + \lambda_1 (x(0)-x_i) + \lambda_2 (x(1)-x_f).\label{eq:BVaction}
\end{align}
We have added two Lagrange multipliers that are treated as additional dynamical degrees of freedom of our system and in turn enforce the boundary conditions of the numerical solution. This procedure may appear to introduce the boundary conditions strongly, however it does not amount to an apriori replacement of $x(t_1)$ and $x(t_2)$ by $x_i$ and $x_f$. During the procedure to locate the critical point of $\mathds{S}_{\rm BVP}$, we find that the minimization algorithms approach the extremum of the functional globally. I.e. the boundary conditions are fulfilled to machine precision for the actual solution, while deviations are possible at intermediate steps. 

Note that when adding Lagrange multipliers to introduce constraints to an optimization functional, the relevant extremum may become a saddle point. If only algorithms are available that locate the minima of a functional, one can circumvent this issue by optimizing the norm of the gradient of the functional instead, for which the saddle point becomes a local minimum. In this study we use as a preconditioning step the gradient-free Nelder-Mead minimizer with a subsequent application of the Newton method and the Interior Point method \footnote{Standard implementations of the aforementioned algorithms in the \texttt{NMinimize} and \texttt{FindMinimum} functions provided by the software \texttt{Mathematica 12.0}  \cite{Mathematica} have been used.}.

\begin{figure}
\includegraphics[scale=0.20]{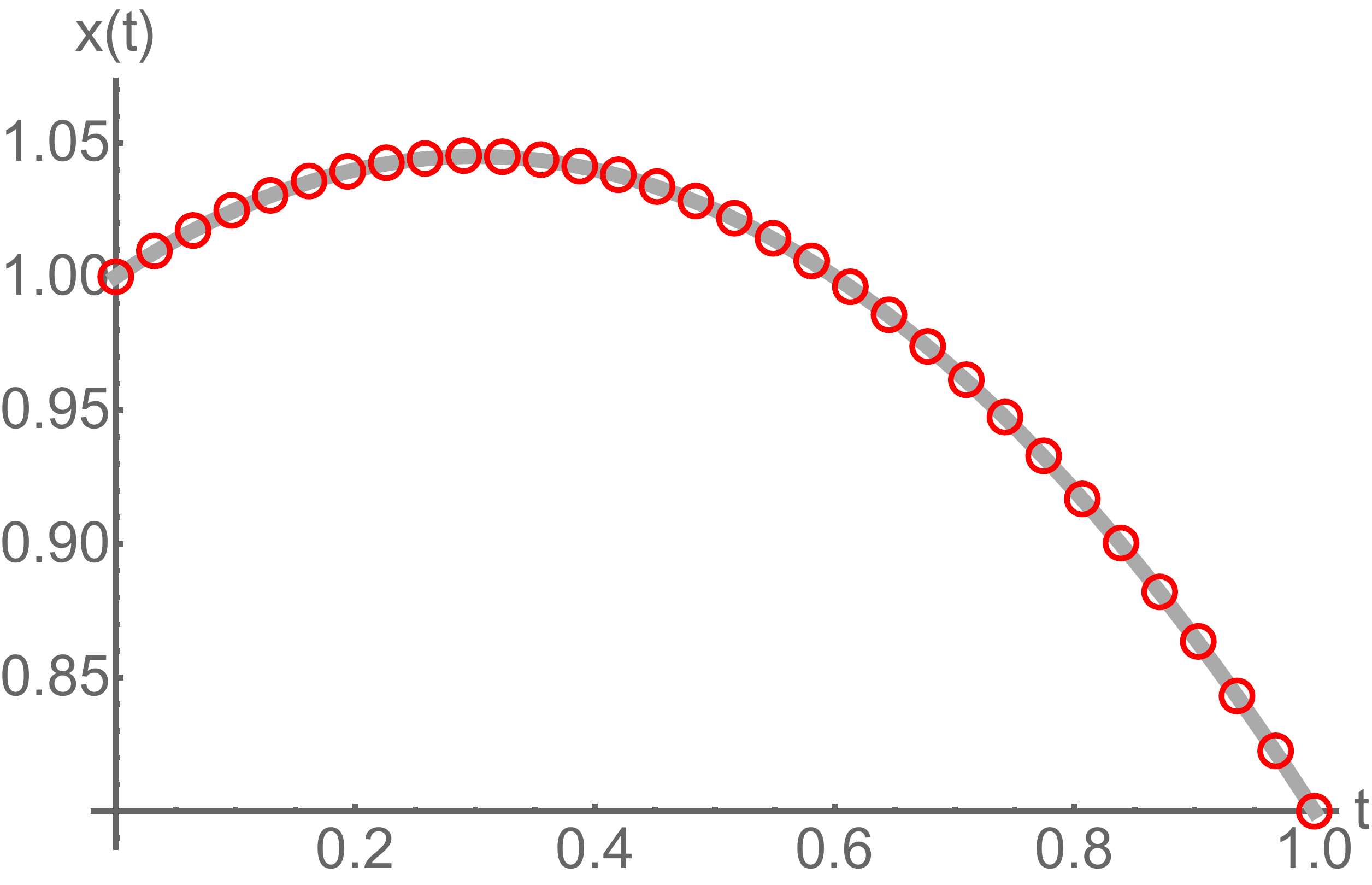}
\includegraphics[scale=0.17]{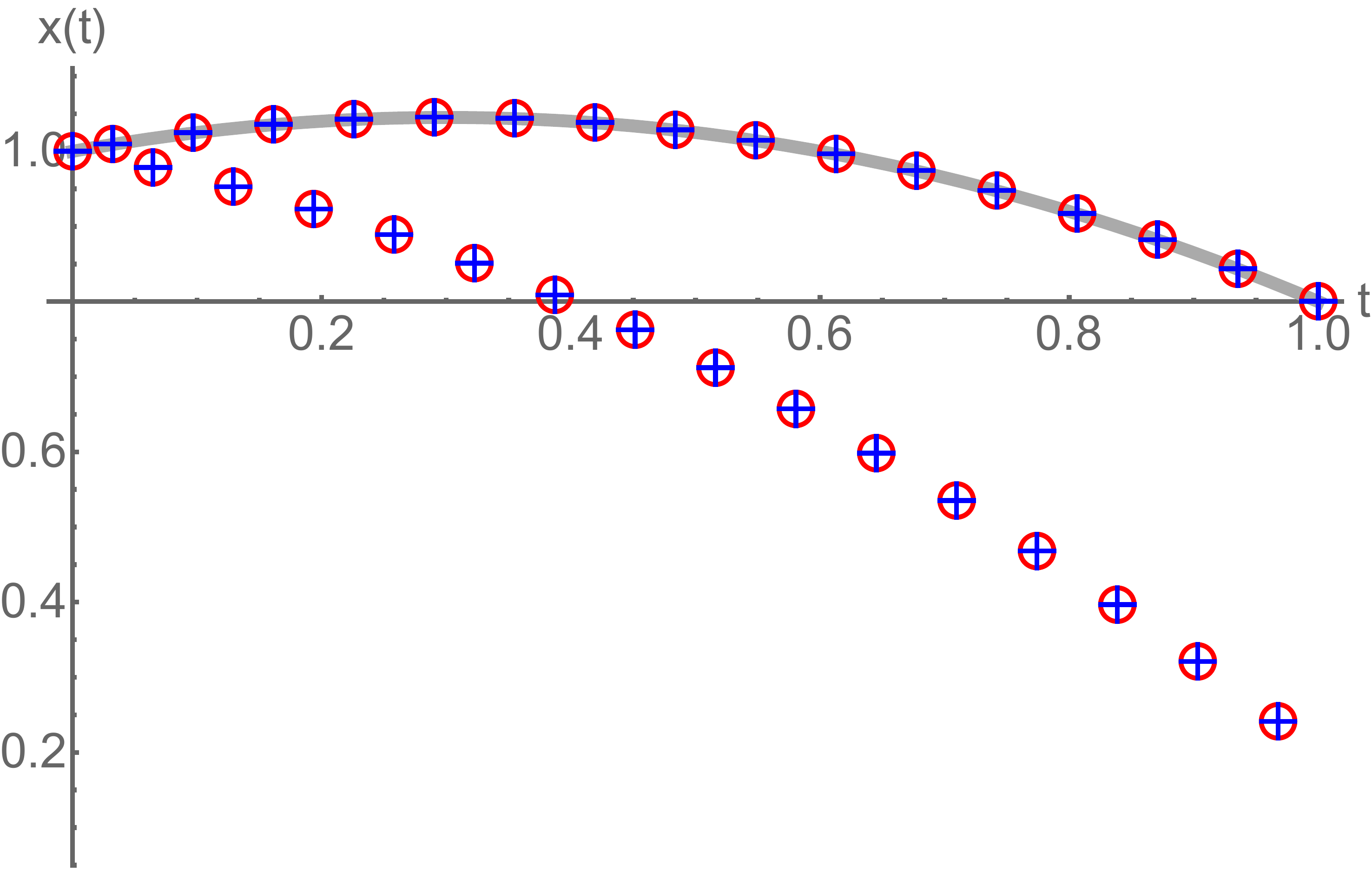}
    \caption{ Using a \texttt{SBP21} operator in time for the boundary value formulation we find (left) the discretized path ${\bf x}$ (red circles) ($N_t=32$) that optimize the functional \cref{eq:BVaction}. The discretized path ${\bf x}_1$ (red circles) and ${\bf x}_2$ (blue crosses) ($N_t=32$) that optimize the functional \cref{eq:IVPfunc} for the initial value formulation, corresponding to the discretized $S_{\rm IVP}$ is shown in the right panel. Note that only half of the path elements reproduce the correct solution. Continuum solution of the Euler-Lagrange equation $\ddot x(t)=-g$ is shown as solid gray line.}
    \label{fig:classtrajBVP}
\end{figure}

Having derived the continuum Euler-Lagrange equations for the point mass in a constant gravitational field before in \cref{eq:ELpointmass}, we compute the explicit solution for the classical trajectory in the time interval $t\in[0,1]$, based on a value of $m=1$, $g=1$ and initial conditions $x(0)=1$, $\dot x(0)=0.3$. In that case the point mass reaches the position $x(1)=0.8$ at time $t_2=1$. Supplying these values to \cref{eq:BVaction}, we can solve for the extremum and, as shown the left panel of \cref{fig:classtrajBVP}, obtain a solution (red dots) that recovers the correct solution of the Euler-Lagrange equations (gray solid line).

While we succeed in recovering the correct solution, this approach, as mentioned before, is conceptually not satisfactory, since the formulation of the BVP relied on information about $x(1)$ obtained from the prior solution of the Euler-Lagrange equations as initial value problem.

\subsection{A naive SBP discretization of the model initial value problem}

Let us continue by turning our attention to discretizing the continuum formulation of the variational principle for initial value problems, which is based on two paths. Introducing discretized paths  ${\bf x}_1=(x_1(0),x_1(\Delta t),x_1(2\Delta t),\ldots)^{\rm T}$ and correspondingly ${\bf x}_2$ and using the same symbols as before for the integration $\mathds{H}$ and summation-by-parts difference operators $\mathds{D}$, we arrive at the following action
\begin{align}
\nonumber \mathds{S}_{\rm IVP}&= \Big\{  \frac{1}{2} (\mathds{D}{\bf x}_1)^{\rm T} \mathds{H} (\mathds{D}{\bf x}_1) - g \mathds{1}^{\rm T} \mathds{H} {\bf x}_1\Big\} - \Big\{\frac{1}{2} (\mathds{D}{\bf x}_2)^{\rm T} \mathds{H} (\mathds{D}{\bf x}_2) - g \mathds{1}^{\rm T}  \mathds{H} {\bf x}_2 \Big\}\\
\nonumber &+ \lambda_1 (x_1(0)-x_i) + \lambda_2((\mathds{D}{\bf x}_1)(0)-\dot x_i) \\
&+ \lambda_3 (x_1(N_t)-x_2(N_t)) + \lambda_4 ( (\mathds{D}{\bf x}_1)(N_t)-(\mathds{D}{\bf x}_2)(N_t) ).\label{eq:IVPfunc}
\end{align}
Here we have introduced four Lagrange multipliers to both enforce the initial conditions for position $x_i$ and derivative $\dot x_i$ of the forward path $(\lambda_1,\lambda_2)$, as well as to enforce the correct identification of the position and derivatives at the last point of the forward and backward path $(\lambda_3,\lambda_4)$.  We consider all $\lambda_i$'s as dynamical degrees of freedom, such that the constraints are enforced exactly on the final solution of the optimization problem, while permitting deviations from the constraints at intermediate steps.

Before we continue to determine the optimal paths according to $\mathds{S}_{\rm IVP}$, we show that this discrete functional yields the correct equations of motion according to the stationarity condition \cref{eq:varIVPdefEq}, if SBP operators are used. Focusing on the terms in curly brackets in \cref{eq:IVPfunc}, we introduce the discretized ${\bf x}_1={\bf x}_+ +(1/2) {\bf x}_-$ and ${\bf x}_2={\bf x}_+-(1/2){\bf x}_-$ so that
\begin{align}
\nonumber \mathds{S}&= \frac{1}{2} (\mathds{D}{\bf x}_+)^{\rm T} \mathds{H} (\mathds{D}{\bf x}_-) + \frac{1}{2} (\mathds{D}{\bf x}_-)^{\rm T} \mathds{H} (\mathds{D}{\bf x}_+) - g \mathds{1}^{\rm T}  \mathds{H} {\bf x}_-\\
\nonumber &=  (\mathds{D}{\bf x}_+)^{\rm T} \mathds{H} (\mathds{D}{\bf x}_-) - g \mathds{1}^{\rm T}  \mathds{H} {\bf x}_-\\
&= -(\mathds{D}\mathds{D}{\bf x}_+)^{\rm T} \mathds{H} {\bf x}_- + (\mathds{D}{\bf x}_+)^{\rm T} (\mathds{E}_N-\mathds{E}_0) {\bf x}_- - g \mathds{1}^{\rm T}  \mathds{H} {\bf x}_-\label{eq:discrEOMderiv}
\end{align}  
Here we have used the symmetry of $\mathds{H}$ to arrive at the second line and explicitly exploited the SBP property of $\mathds{D}$ in the third line. As we enforce the initial conditions and identify the forward and backward path at the final time, both boundary terms involving ${\bf x}_- $ vanish. Mimicking the continuous derivation, let us take the derivative of $\mathds{S}$ with respect to the i-th component of the vector ${\bf x}_-$, which yields the following expression
\begin{align}
\Big( -(\mathds{D}\mathds{D}{\bf x}_+)^{\rm T} - g \mathds{1}^{\rm T} \Big)\mathds{H}{\bf e}_i =0. \label{eq:discreom2}
\end{align}
Since $\mathds{H}$ is diagonal, \cref{eq:discreom2} establishes the discrete equation of motion $(\mathds{D}\mathds{D}{\bf x}_+) = - g \mathds{1}$, a faithful representation of the continuum result $\ddot x(t)=-g$.

Let us continue to determining the optimal paths ${\bf x}_1$ and ${\bf x}_2$ that correspond to the critical point of $\mathds{S}_{\rm IVP}$ using $\mathds{H}^{[2,1]}$ and $\mathds{D}^{[2,1]}$ on $N_t=32$ discrete points, we find the solution shown in the bottom panel of \cref{fig:classtrajBVP}. We plot the values of the forward path as red circles, while those of the backward path are given as blue crosses. Note that they lie on top of each other, which tells us that the optimal solution fulfills the \textit{physical limit} condition ${\bf x}_1={\bf x}_2$.

On the other hand we also immediately see that only around half of the points on each path agree with the correct solution from the Euler-Lagrange equations (gray solid line). The other half lies significantly below the correct solution, forming a highly oscillatory structure. For an even number of grid points the last point of ${\bf x}_1$ and ${\bf x}_2$ lies on the correct trajectory, while for an odd number of points, the path ends on the oscillatory structure below.

We have identified the origin of these oscillatory structures to arise from the particular structure of the null-space of the finite difference operator. In the kinetic terms of $\mathds{S}_{\rm IVP}$ both $\mathds{D}$ and $\mathds{D}^{\rm T}$ appear. The study of null-space consistency of the lowest order \texttt{SBP21} operator considered here, reveals that it contains exactly two zero eigenvalues. The space of right eigenvectors of $\mathds{D}$, associated with this doubly degenerate eigenvalue, is only one-dimensional. Both eigenvectors are proportional to the constant function. 

\begin{figure}
    \centering
    \includegraphics[scale=0.28]{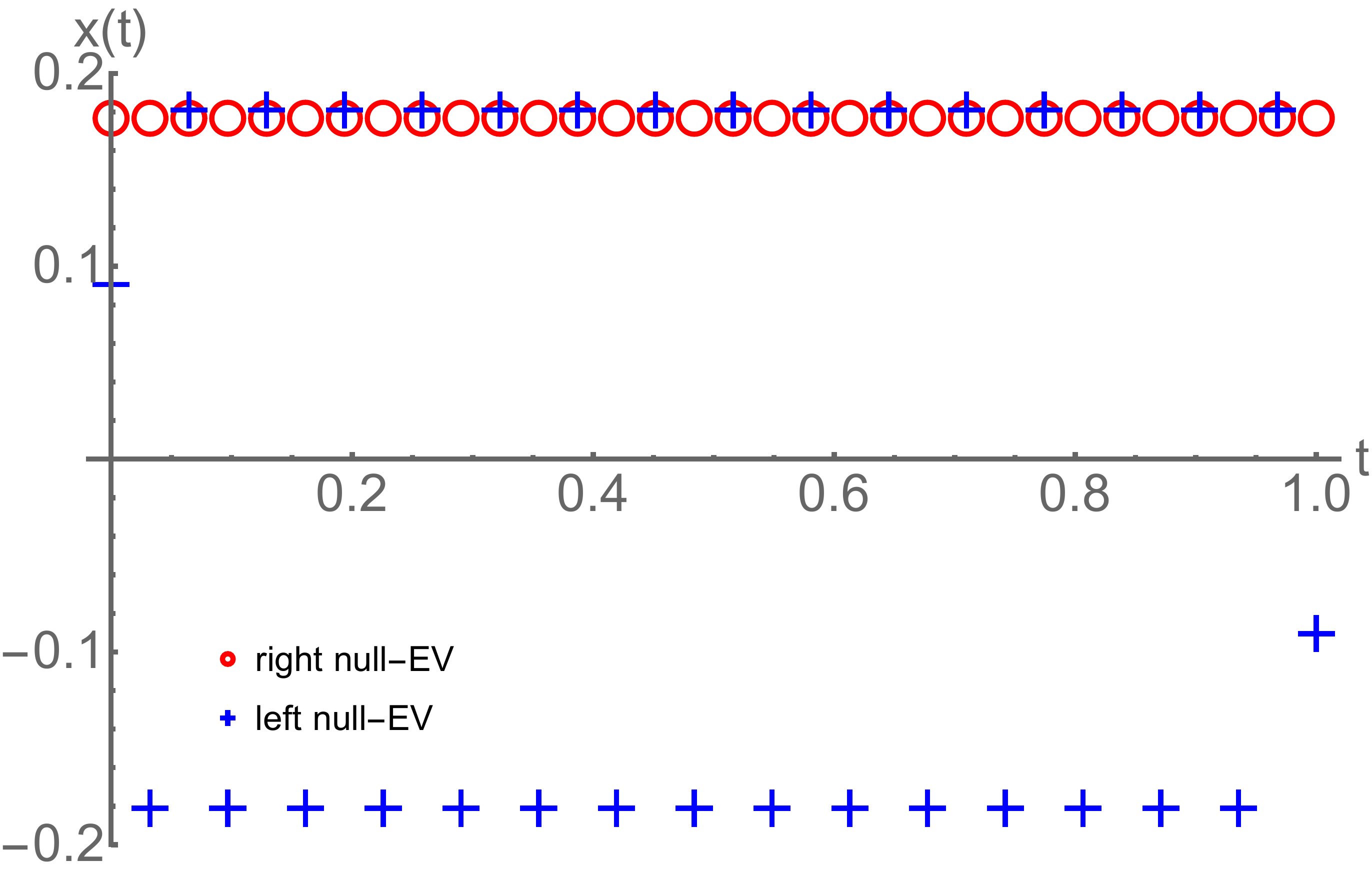}
    \caption{Right (red circles) and left (blue crosses) eigenvector associated with a zero eigenvalue in $\mathds{D}^{[2,1]}$, based on $N_t=32$ points. Note the highly oscillatory character of the latter.}
    \label{fig:DEigenvec}
\end{figure}

However when we study the form of the left eigenvectors of $\mathds{D}$, or equivalently the right eigenvectors of $\mathds{D}^{\rm T}$, we find that those projecting into the null space are not at all constant but highly oscillatory, reminiscent of the so-called $\pi$-mode. An example of these eigenvectors is shown in \cref{fig:DEigenvec}.

Such unphysical oscillatory solutions have recently been identified to also interfere in determining the solutions of differential equations with non-trivial boundary conditions in one- and multiple dimensions in \cite{ranocha_discrete_2020}. In the context of the variational problem considered here, the oscillatory solutions did not affect the solution when formulated as a boundary value problem in \cref{sec:discbvp}. The fixing of the boundary at $t_2$ in the BVP formulation apparently prevents the oscillatory solution. On the other hand the IVP action \cref{eq:IVPfunc} clearly accommodates these oscillatory paths.

An accurate discretization scheme for the IVP system action must therefore be able to avoid the appearance of unphysical oscillatory modes and several strategies to do so have been explored in the literature. One class of strategies consists of modifying the first order derivative operator $\mathds{D}$ by adding higher order derivative operators to it. A conventional SBP finite difference operator of order $p$ needs to fulfill the derivative property $\mathds{D}^{(p)} {\bf x}^n = (n-1){\bf x}^{n-1}$ exactly only on monomials up to order $n=p-1$. Therefore, adding a higher derivative operator $\Delta t \, \mathds{D}^{(p+2)}$, scaled by the grid spacing does not affect this property, as it annihilates all lower order monomials. In addition, this correction term vanishes in the limit of taking $\Delta t\to0$. In the context of upwind schemes one e.g. adds the  symmetric second derivative operator to the SBP first derivative, turning it into an upwind derivative. If one deals with complex functions one may instead add the symmetric second derivative multiplied with the imaginary unit. This modification is known as adding a \textit{Wilson term} \cite{wilson_confinement_1974} in the physics literature\footnote{Wilson derived that regularization after investigating the Green's function of the differential operator that defines the equation of motion of the system. What he found is that in Fourier space the Green's function exhibits not only a pole corresponding to the physical trajectory, but due to the finite grid spacing a second pole appeared at the end of the Brillouin zone, which introduces exactly the oscillatory mode we observed above.}.

Both of these approaches present challenges, which we wish to avoid here. By turning the central stencil into an upwind stencil, we lose the symmetry of the system, which adversely affects the accuracy of the solutions. Introducing a purely imaginary modification on the other hand requires the difference operator to act on complex functions to be meaningful. One may contemplate the possibility to complexify the functions involved in the variational problem, which while only cumbersome in the classical case will lead to conceptual problems when trying to use the discretization in the context of quantum path integrals (c.f. sign problem).

We therefore wish to explore a different route to remove the unphysical zero modes of the operator $\mathds{D}$, taking inspiration from more recent works on null-space consistent SBP operators, such as in refs.~\cite{svard2019convergence,linders_properties_2020,svard_convergence_2021,ranocha2021new}. The central ingredient in these approaches is to exploit the weak formulation of boundary and initial conditions. Concretely, when boundary conditions are enforced weakly via a penalty term, this penalty term can be partially absorbed into the derivative operator to remove the zero modes of that operator. On the level of differential equations, the strategy works as follows. Consider the following IVP, the differential equation for the exponential function
\begin{align}
\frac{d}{dt} u(x)=\lambda u(x),\quad u(0)=u_0,
\end{align}
which in its discretized form reads 
\begin{align}
\mathds{D} {\bf u} = \lambda {\bf u} + \sigma_0 \mathds{H}^{-1} \mathds{E}_0\big( {\bf u} - {\bf g}\big). \label{eq:simpleODE}
\end{align}
Here we have added a so-called SAT penalty term on the RHS, which includes the matrix $\mathds{E}_0={\rm diag}[1,0,\ldots,0]$ that singles out the first entry in the discretized functions ${\bf u}$ and ${\bf g}$. The former ${\bf u}$ refers to the solution of the differential equation and the latter ${\bf g}=(u_0,0,\ldots,0)$ contains the initial value as its first entry. Note that $\mathds{H}^{-1}$ contains ${\Delta t}^{-1}$, which contributes with increasing weight as ${\Delta t}\to0$. The parameter $\sigma_0$ in the SBP-SAT approach is tuned to satisfy stability properties and its optimal value is found to be $\sigma_0=-1$, a choice we adopt in the following. The standard approach developed in the conventional SBP-SAT treatment of IVPs consists of absorbing the penalty term proportional to ${\bf u}$ into a redefined $\tilde{\mathds{D}}=\mathds{D} - \sigma_0 \mathds{H}^{-1} E_0$, which does not feature any zero modes anymore. That operator is now non-singular \cite{ruggiu_eigenvalue_2020} and may be inverted to obtain the solution ${\bf u}$. In the next section we will develop a similar strategy applicable to the variational problem.

\subsection{Regularized SBP discretization of the model initial value problem}
\label{sec:strongIVP}

Taking inspiration from the work on regularizing SBP operators in differential equations, we set out to absorb information about the initial conditions into the SBP operator as means of regularization. In the functional of \cref{eq:IVPfunc} we do not have an equality sign, such as in our example \eqref{eq:simpleODE}, to rearrange terms. Instead we must find a way to incorporate the whole penalty term in $\mathds{D}$. Note that the penalty term contains one expression that is proportional to the function that the SBP operator acts on and one expression proportional to a constant. I.e. we have to modify the difference operator to include a shift. In other words, we are dealing with an affine transformation.

\begin{figure}
    \includegraphics[scale=0.3]{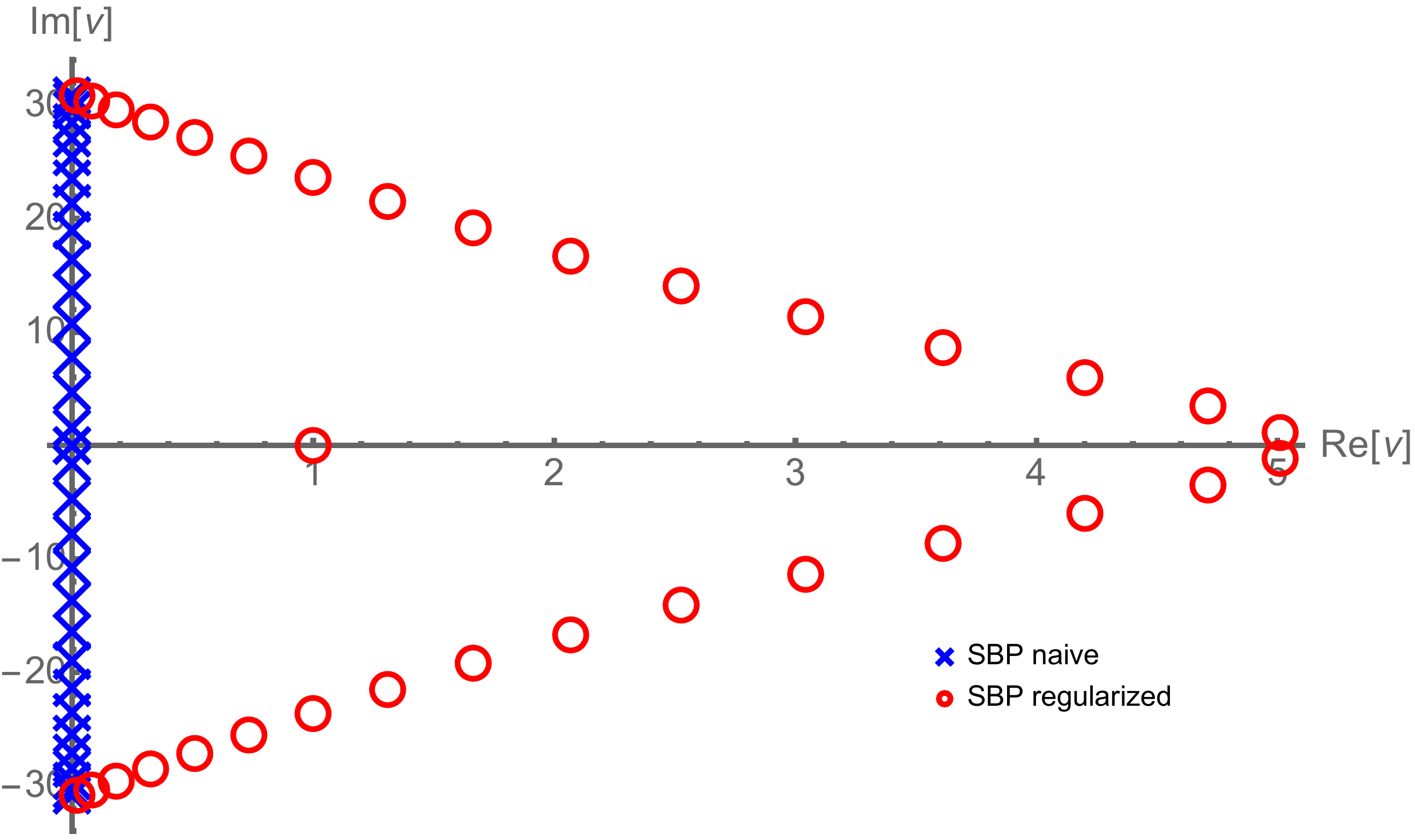}
    \includegraphics[scale=0.33]{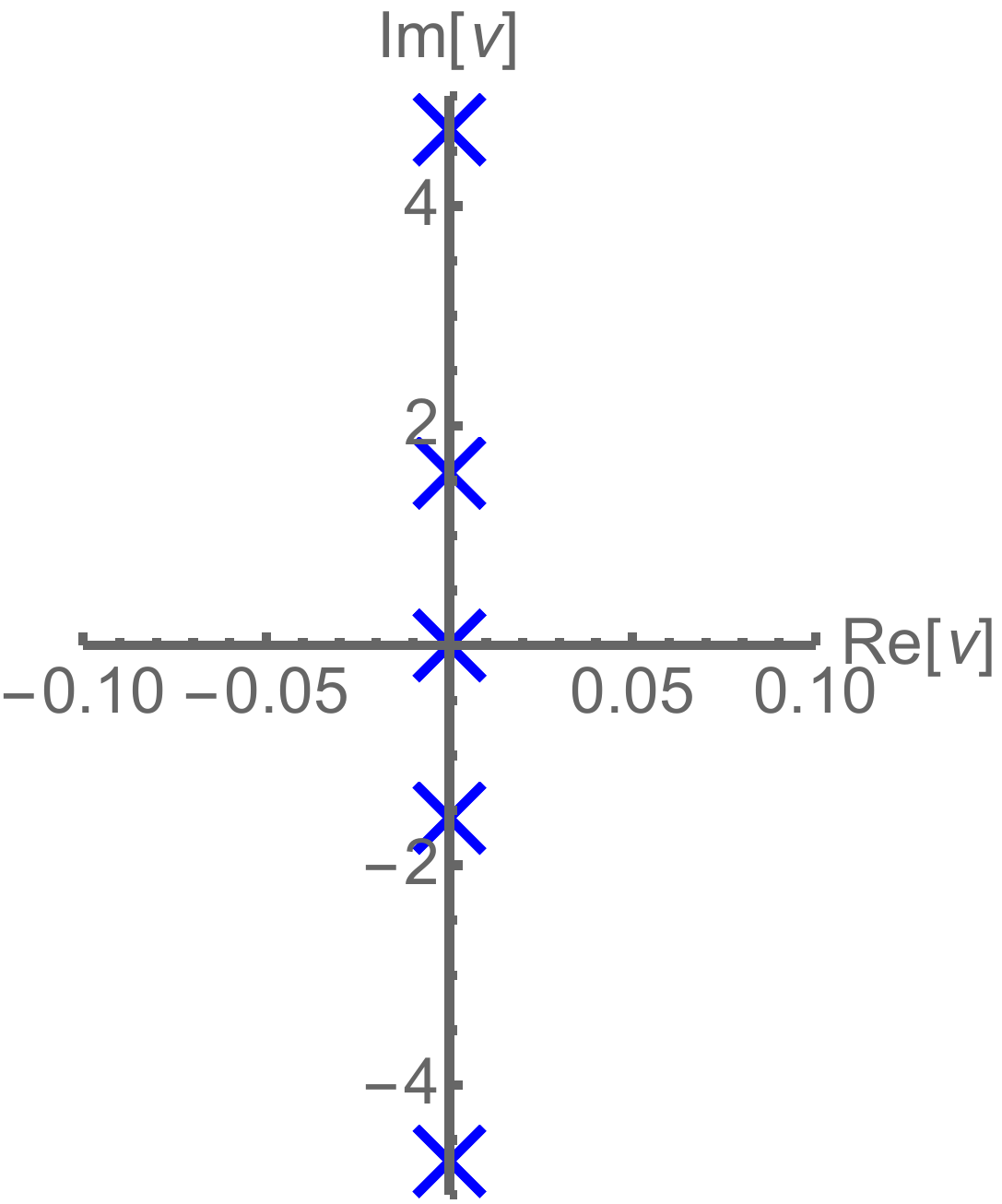}
    \caption{(Left) Eigenvalue spectrum of the unregularized \texttt{SBP21} operator $\mathds{D}^{[2,1]}$ (blue crosses) for $N_t=32$ grid points and corresponding spectrum of the regularized operator $\bar{\mathds{D}}^{[2,1]}$ (red circles), which does not feature any zero modes. The zoomed inset on the right reveals the presence of the zero modes in the unregularized \texttt{SBP21} operator $\mathds{D}^{[2,1]}$. }
    \label{fig:DEigenvalsmod}
\end{figure}

There exists an elegant way to express affine transformations using so-called affine coordinates. One defines $\bar A[{\bf b}] \bar{\bf x} = A{\bf x}+{\bf b}$, where $\bar A[{\bf b}]$ refers to the matrix $A$ amended by one more row and column with $1$ placed in the lower right corner. The additional column available in $\bar A[{\bf b}]$ is filled with the values of ${\bf b}$. The vector $\bar{\bf x}$ is just ${\bf x}$ amended by one more entry with value one. For our application to the variational formulation of the IVP we therefore define a new $\bar{\mathds{D}}$ using as shift the vector containing the initial values ${\bf b}= \sigma_0 \mathds{H}^{-1} E_0 {\bf g}$ where ${\bf g}={\rm diag}[x_i,x_i+\Delta t\, \dot x_i,0,\cdots,0]$. For the \texttt{SBP21} operator\footnote{For a higher order SBP operator, the values of ${\bf g}$ need to be chosen, so that $[\mathds{D} {\bf g}](0)=\dot x_i$ and $[{\bf g}](0)=x_i$.} the explicit expression we obtain reads
\begin{align}
\bar{\mathds{D}}^{[2,1]}=
\left[ \begin{array}{cccccc} -\frac{1}{\Delta t} -\sigma_0 \frac{2}{\Delta t}&\frac{1}{\Delta t} & & &&\sigma_0 \frac{2}{\Delta t} x_i\\ -\frac{1}{2\Delta t}& 0& \frac{1}{2\Delta t}& &&0\\ & &\ddots &&&\vdots\\ &&-\frac{1}{2\Delta t}&0&\frac{1}{2\Delta t}&0\\ &&&-\frac{1}{\Delta t}&\frac{1}{\Delta t}&0\\  0&&\ldots &&0&1\\  \end{array} \right].
\end{align}
Note that in this paper we choose the parameter $\sigma_0=-1$, whenever a penalty term is incorporated in ${\bar {\mathds{D}}}$. This choice is motivated by the fact that in the conventional treatment of IVPs using the SBP-SAT approach, this value leads to a minimal discretization error (see e.g. ref.~\cite{lundquist2014sbp}). We find that $\sigma_0=-1$ is optimal for our approach too, as only in this case the correct classical solution is recovered.

All zero modes of the original operator $\mathds{D}^{[2,1]}$ have been lifted in $\bar{\mathds{D}}^{[2,1]}$ and the resulting spectrum of eigenvalues $\nu$ is shown in \cref{fig:DEigenvalsmod}. Note that $\bar{\mathds{D}}^{[2,1]}$ still correctly annihilates the constant function, as long as it is compatible with the initial conditions $x(0)=x_i$. In affine coordinates this annihilation does not lead to a resulting zero vector, but a vector that contains vanishing entries, except for the final one associated with the single real eigenvalue of value one, shown in \cref{fig:DEigenvalsmod}.

\begin{figure}
    \centering
    \includegraphics[scale=0.28]{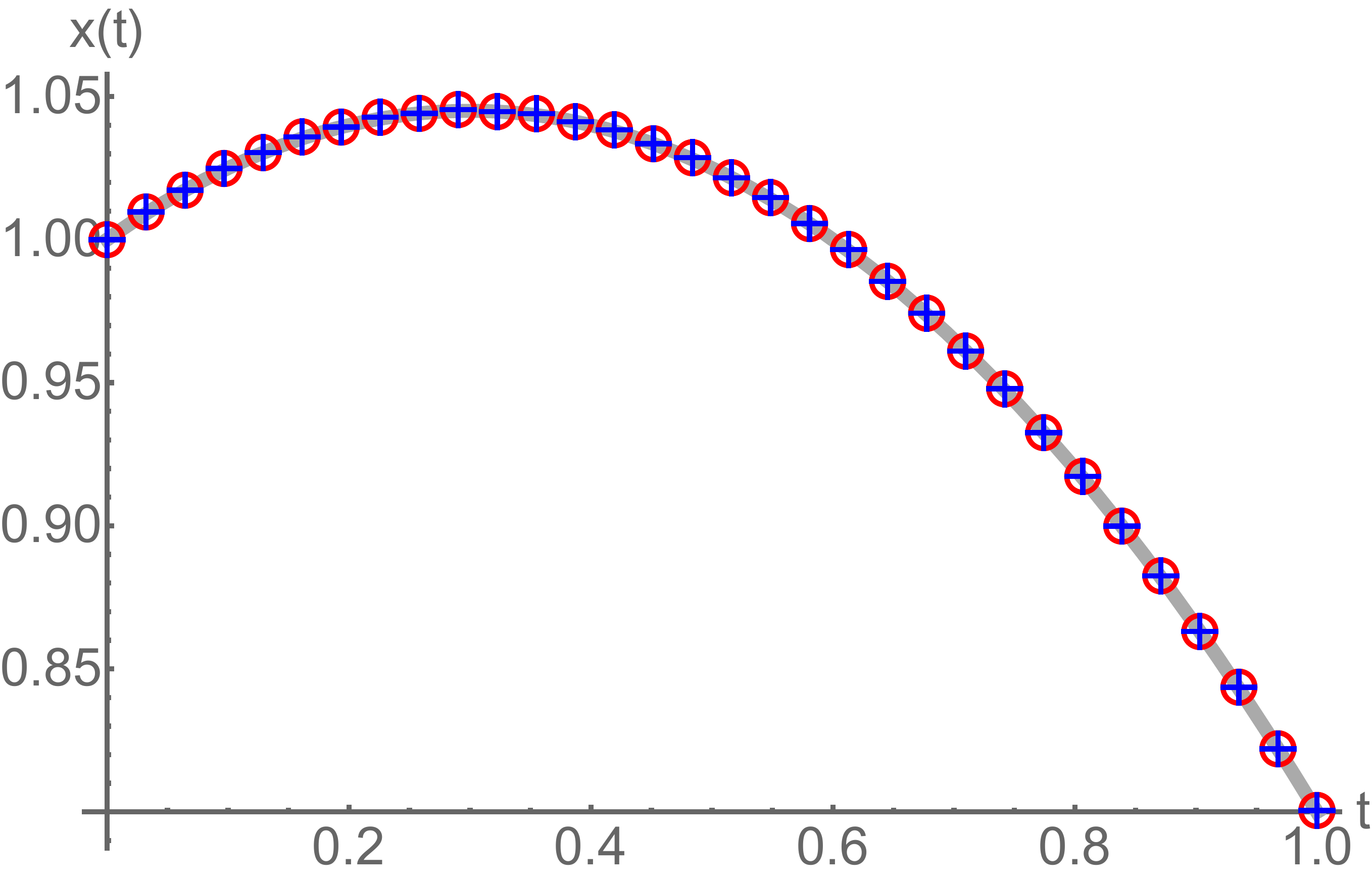}
    \caption{Numerical solution of the discretized path ${\bf x}_1$ (red circles) and ${\bf x}_2$ (blue crosses) ($N_t=32$) that optimize the functional \cref{eq:IVPfuncReg} corresponding to the discretized $S_{\rm IVP}$. We use the regularized \texttt{SBP21} operator in time. Continuum solution of the Euler-Lagrange equation $\ddot x(t)=-g$ is shown as solid gray line. Note that the solution successfully avoids oscillatory contamination.}
    \label{fig:classtrajIVPSBPReg}
\end{figure}
When formulating the action with the modified SBP operator, we obtain
\begin{align}
\nonumber \mathds{S}_{\rm IVP}&= \Big\{  \frac{1}{2} (\bar{\mathds{D}}{\bar{\bf x}}_1)^{\rm T} \bar{\mathds{H}} (\bar{\mathds{D}}\bar{{\bf x}}_1) - g \mathds{1}^{\rm T} {\mathds{H}} {\bf x}_1 \Big\} - \Big\{\frac{1}{2} (\bar{\mathds{D}}\bar{{\bf x}}_2)^{\rm T} \bar{\mathds{H}} (\bar{\mathds{D}}\bar{{\bf x}}_2) - g \mathds{1}^{\rm T} \mathds{H} {\bf x}_2 \Big\}\\
\nonumber &+ \lambda_1 (x_1(0)-x_i) + \lambda_2((\mathds{D}{\bf x}_1)(0)-\dot x_i) \\
&+ \lambda_3 (x_1(N_t)-x_2(N_t)) + \lambda_4 ( (\mathds{D}{\bf x}_1)(N_t)-(\mathds{D}{\bf x}_2)(N_t) ).\label{eq:IVPfuncReg}
\end{align}
In order to implement the inner product in affine coordinates, we define $\bar{\mathds{H}}$, which denotes the matrix $\mathds{H}$, amended by one extra row and column of values zero. The last entry of the vector $\bar{\mathds{D}}{\bar{\bf x}}_{1,2}$ serves only to implement the shift in affine coordinates, hence it can be discarded via $\bar{\mathds{H}}$ since the regularized SBP operator has already acted on the path.

Note that here we again add Lagrange multipliers as dynamical degrees of freedom, to fulfill the initial conditions. One may ask whether enforcing the initial conditions in this way neutralizes the effect of the regularization. We emphasize that this is not the case. Minimization algorithms approach the extremum of the functional globally, allowing the regulator to remain effective and to avoid the oscillatory solutions.

Another question of both conceptual and practical relevance is whether the functional \cref{eq:IVPfuncReg} houses one or multiple different local extrema. For the case of the point particle in a constant gravitational field, we find that the answer is that the functional is convex and a thus any local extremum is also a global extremum. Let us determine the curvature of \cref{eq:IVPfuncReg} with respect to the individual entries of the paths ${\bf x}_{1,2}$. The matrix $A_{ij}=\partial^2\mathds{S}_{\rm IVP} /\partial x_i \partial x_j=[ \bar{\mathds{D}}^{\rm T} \bar{\mathds{H}}\bar{\mathds{D}}]_{ij}$ is indeed positive semi-definite, as can be checked explicitly using a computer algebra system. However, in general convexity is not automatic and needs to be checked on a case-by-case basis\footnote{None of the functionals treated in this study suffered from multiple extrema, allowing the Newton and Quasi-Newton methods implemented in Mathematica \cite{Mathematica} to arrive at a single solution independent of starting point.}. 

For the point particle in a constant gravitational field, the solutions ${\bf x}_1$ and ${\bf x}_2$ obtained with \cref{eq:IVPfuncReg} and the regularized \texttt{SBP21} operator are shown in \cref{fig:classtrajIVPSBPReg} as red circles and blue crosses respectively.\footnote{An explicit implementation of all examples discussed in this manuscript can be found as open-access Mathematica script at the Zenodo repository \cite{rothkopf_mathematica_2022}.}. The regularization has successfully removed the contamination by an unphysical oscillatory mode and we are able to reproduce the correct classical solution.

The initial conditions, implemented using Lagrange multipliers and expressed in a variational formulation, have provided us with \textit{a novel discretization prescription} for a wide range of classical systems whose differential equations of motion contain second order derivatives. We emphasize that we did not have to derive the equation of motion to compute the classical trajectory here.

\begin{figure}
  \includegraphics[scale=0.24]{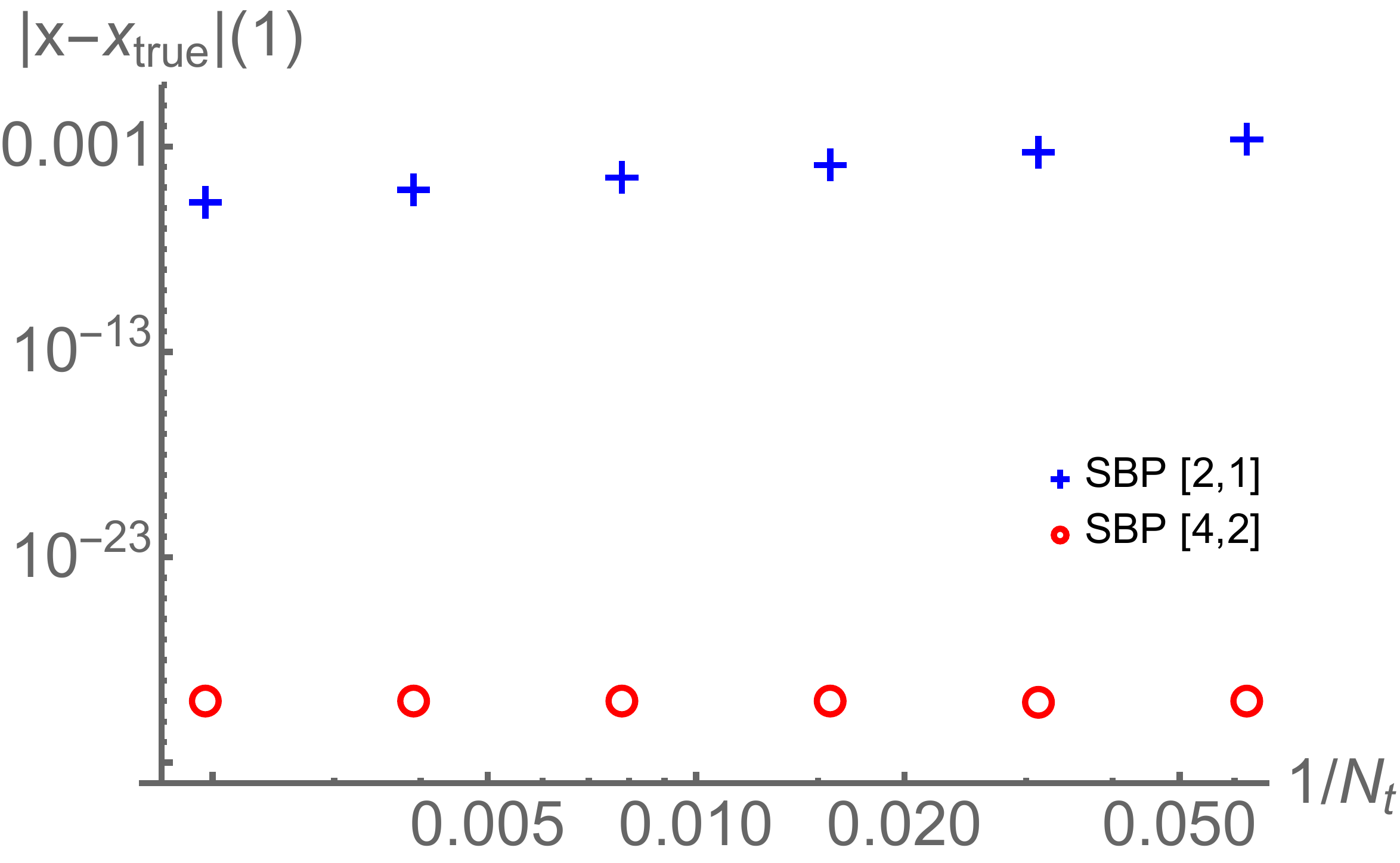}
    \includegraphics[scale=0.24]{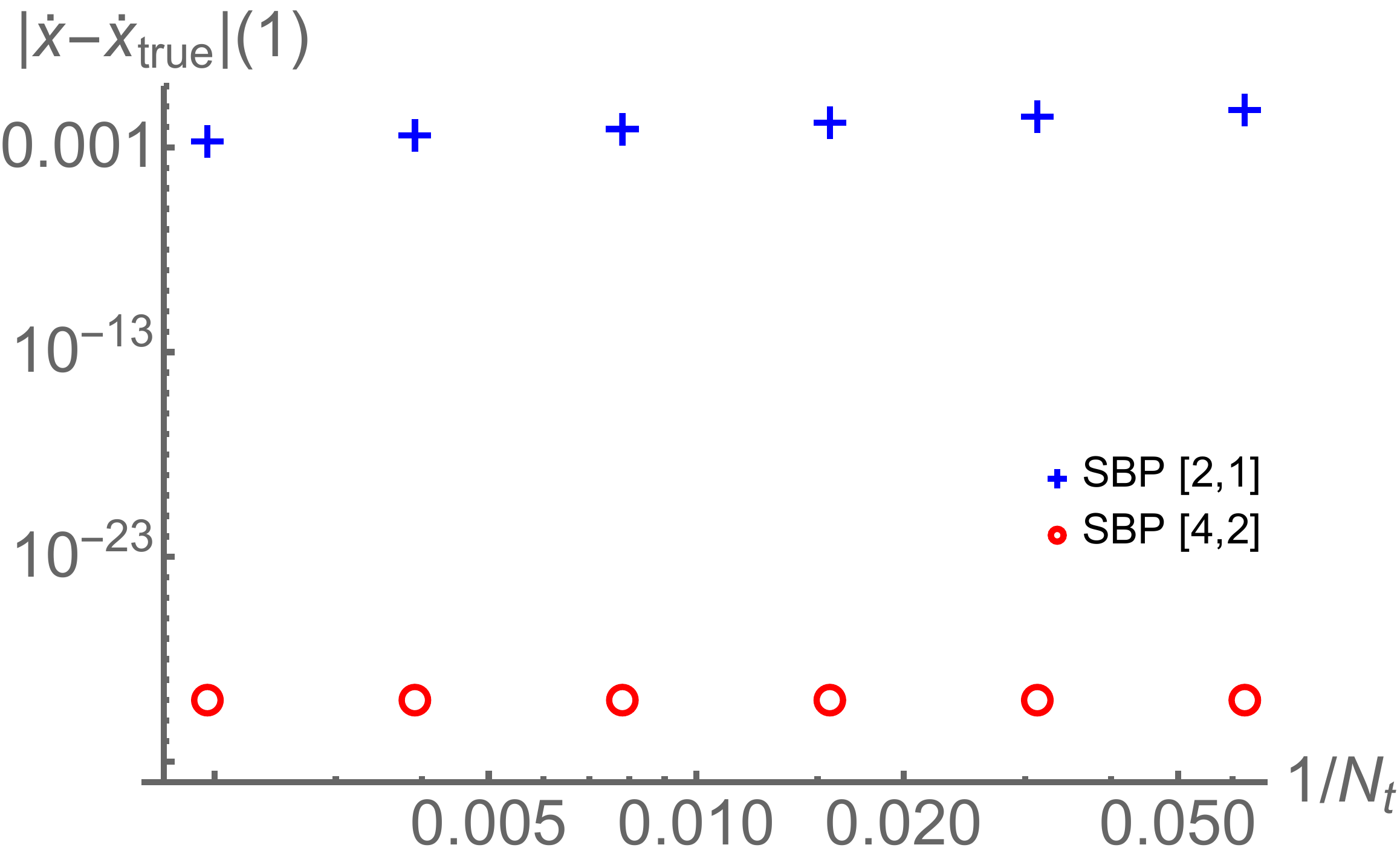}
    \caption{Deviation of the value (left plot) and the derivative (right plot) between the optimal numerical solution of \cref{eq:IVPfuncReg} and the true solution at the final time $t_2=1$. The deviation for the regularized \texttt{SBP21} operator is given as blue crosses, the one for the \texttt{SBP42} operator as red circles. Note that while in the \texttt{SBP21} case the solution improves steadily with diminishing $\Delta t\sim 1/N_t$, the results for the \texttt{SBP42} operator are already exact to machine precision. For a more detailed view of the \texttt{SBP21} behavior see \cref{fig:fitconvsecond}}
    \label{fig:errcompmassgrav}
\end{figure}

Let us take a look at the accuracy and convergence properties of the discretization scheme constructed in this section. The simple model of a point mass in a constant gravitational field again serves as explicit example. To this end we compute the optimal path ${\bf x}_1$ according to the appropriately regularized \cref{eq:IVPfuncReg} using both the regularized \texttt{SBP21} and regularized \texttt{SBP42} operator on different grids using $N_t\in[16,\ldots,512]$ points. We compare the values of the path at the final step $t_2=1$ to the analytically known solution and compute the absolute error between them. These errors are shown in \cref{fig:errcompmassgrav}. 

\begin{figure}
    \centering
    \includegraphics[scale=0.28]{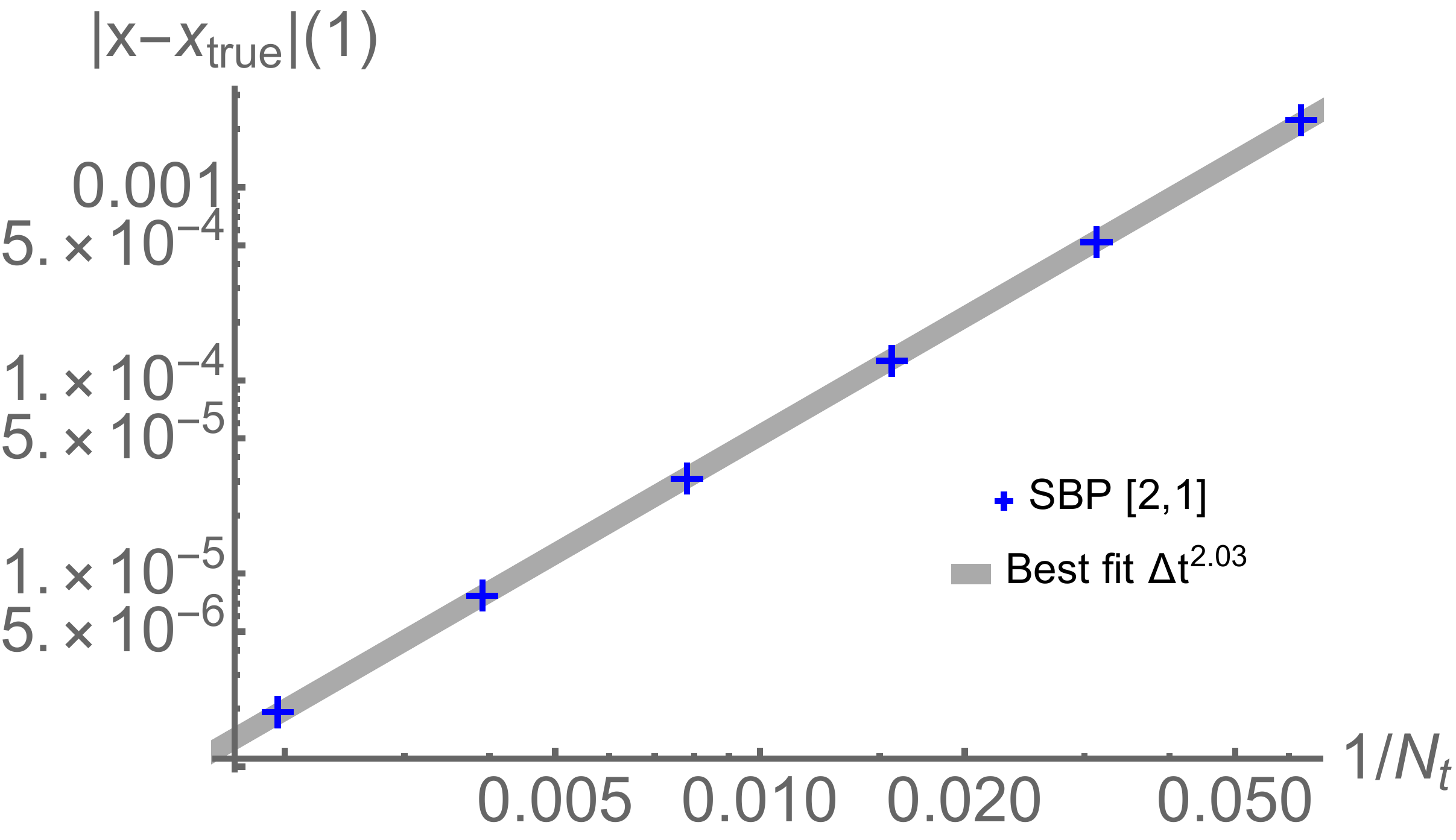}
    \caption{Deviation of the value between the optimal numerical solution of \cref{eq:IVPfuncReg} and the true solution at the final time $t_2=1$ (blue crosses) for the \texttt{SBP21} operator. The best power law fit to the behavior is shown as the gray solid line which corresponds to $\Delta t^{2.03}$.}
    \label{fig:fitconvsecond}
\end{figure}

The lowest order \texttt{SBP21} approximation (blue crosses) exhibits steady improvement in the residual deviation from the true solution, as the grid spacing is reduced. In \cref{fig:fitconvsecond} we zoom in on the \texttt{SBP21} errors and fit them with a power-law ansatz, which on the log-log plot appears as a straight line. The best fit exponent $\Delta t^{2.03}$ we obtain, tells us that our discretization scheme achieves second order accuracy in the solution values. Interestingly, when considering the \texttt{SBP42} operator, we find that irrespective of the grid spacing we are able to reproduce the true solution down to machine precision (which in our case, using the \texttt{Mathematica} software package, was set to $10^{-30}$). This result is reassuring, as by construction the SBP operator and the corresponding quadrature rule are able to differentiate and integrate polynomials up to second order exactly. Since the solution of the point mass in the constant gravitational field is a parabola, we do not find any residual dependence on the grid spacing.

\subsection{Discretization of non-linear second-order IVPs}
\begin{figure}
    \centering
    \includegraphics[scale=0.3]{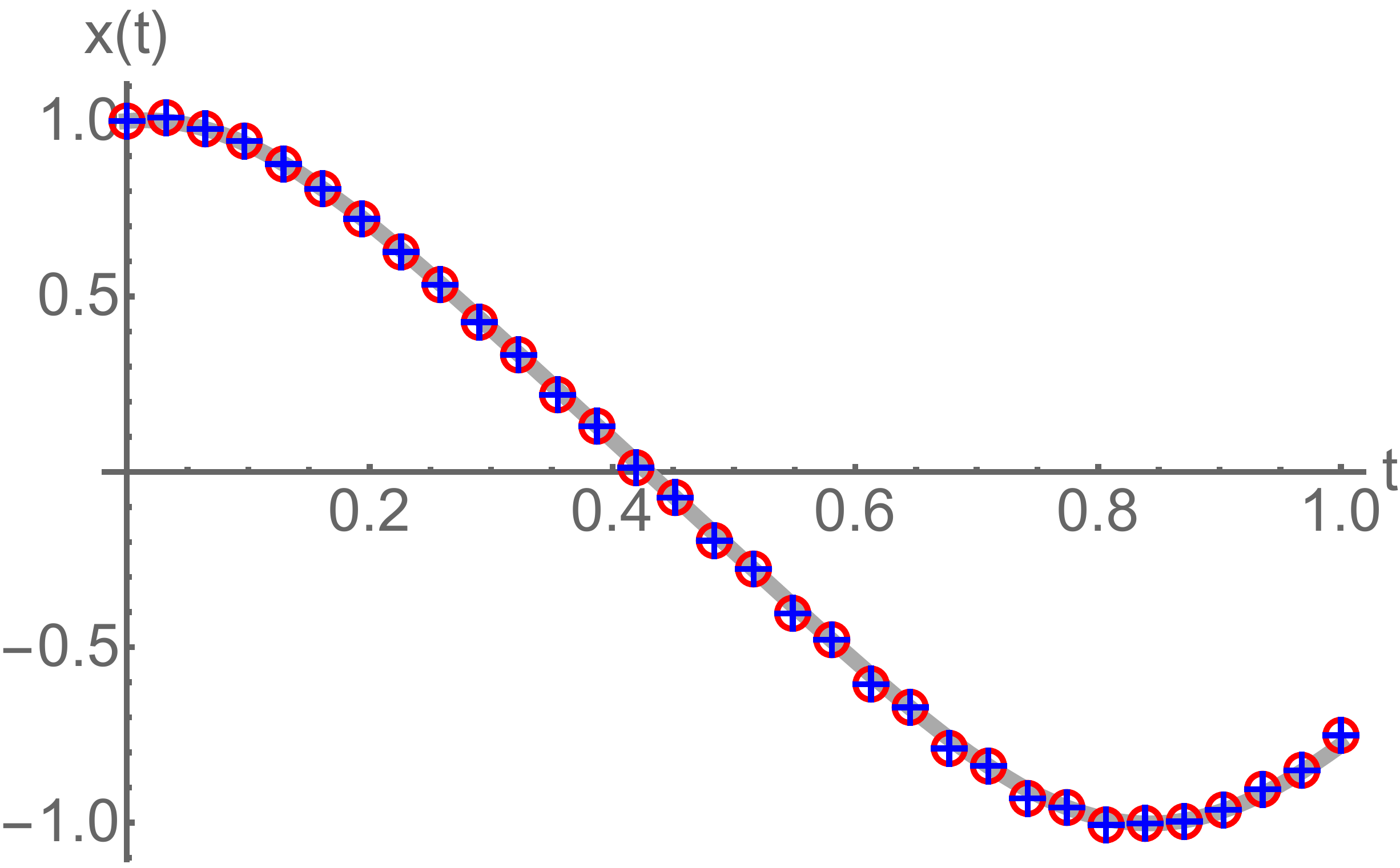}
    \caption{Numerical solution of the discretized path ${\bf x}_1$ (red cricles) and ${\bf x}_2$ (blue crosses) ($N_t=32$) that optimize the functional \cref{eq:IVPnonlinfuncReg} corresponding to the discretized $S_{\rm IVP}$. We use the regularized \texttt{SBP21} summation-by-parts operator in time. Continuum solution of the Euler-Lagrange equation $\ddot x(t) + \kappa x^3(t) =0$ for $\kappa=20$ is shown as solid gray line. Note that the solution successfully avoids oscillatory contamination.}
    \label{fig:classtrajnonlIVPSBPReg}
\end{figure}

While instructive, our model example described a very simple linear system. Let us use the formalism established in \cref{eq:constrsecord} to apply our discretization prescription to a system, which features a genuinely non-linear differential equation of motion instead. The action functional behind the ODE
\begin{align}
\ddot x(t) + \kappa x^3(t) =0,
\end{align} 
is given by the following expression
\begin{align}
{\cal S}_{\rm IVP}=\int dt \Big(\dot x_+(t)\dot x_-(t) - \kappa x_+^3(t)x_-(t)\Big)\label{eq:constrsecord2},
\end{align}
which we must subsequently discretize. The resulting optimization functional reads
\begin{align}
\nonumber \mathds{S}_{\rm IVP}&= \Big\{  \frac{1}{2} (\bar{\mathds{D}}{\bar{\bf x}}_1)^{\rm T} \bar{\mathds{H}} (\bar{\mathds{D}}\bar{{\bf x}}_1)  \Big\} - \Big\{\frac{1}{2} (\bar{\mathds{D}}\bar{{\bf x}}_2)^{\rm T} \bar{\mathds{H}} (\bar{\mathds{D}}\bar{{\bf x}}_2) \Big\}\\
\nonumber & -\kappa\Big(( {{\bf x}}_1 + {{\bf x}}_2)/2\Big)^3\mathds H ( {{\bf x}}_1 - {{\bf x}}_2)\\
\nonumber &+ \lambda_1 (x_1(0)-x_i) + \lambda_2((\mathds{D}{\bf x}_1)(0)-\dot x_i) \\
&+ \lambda_3 (x_1(N_t)-x_2(N_t)) + \lambda_4 ( (\mathds{D}{\bf x}_1)(N_t)-(\mathds{D}{\bf x}_2)(N_t) ).\label{eq:IVPnonlinfuncReg}
\end{align}
The third power in the second line is understood as acting element-wise on the entries of $( {{\bf x}}_1 + {{\bf x}}_2)/2$. The continuum trajectory for the choice $\kappa=20$ is given as the gray solid line in \cref{fig:classtrajnonlIVPSBPReg}. We plot it together with the numerical solutions ${\bf x}_1$ (red circles) and ${\bf x}_2$ (blue crosses) based on the \texttt{SBP21} operator along $N_t=32$ grid points.

\begin{figure}
    \includegraphics[scale=0.24]{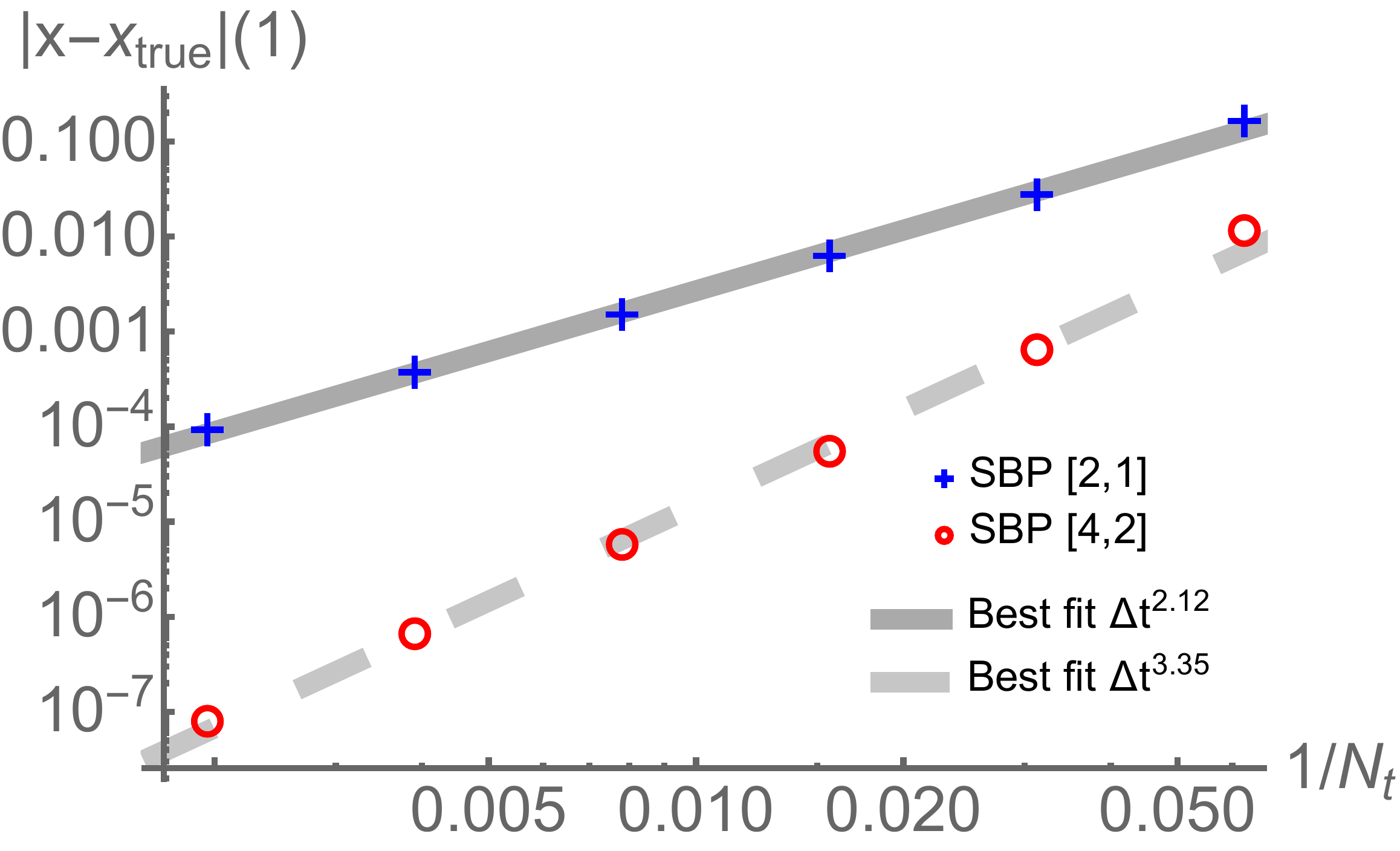}
    \includegraphics[scale=0.24]{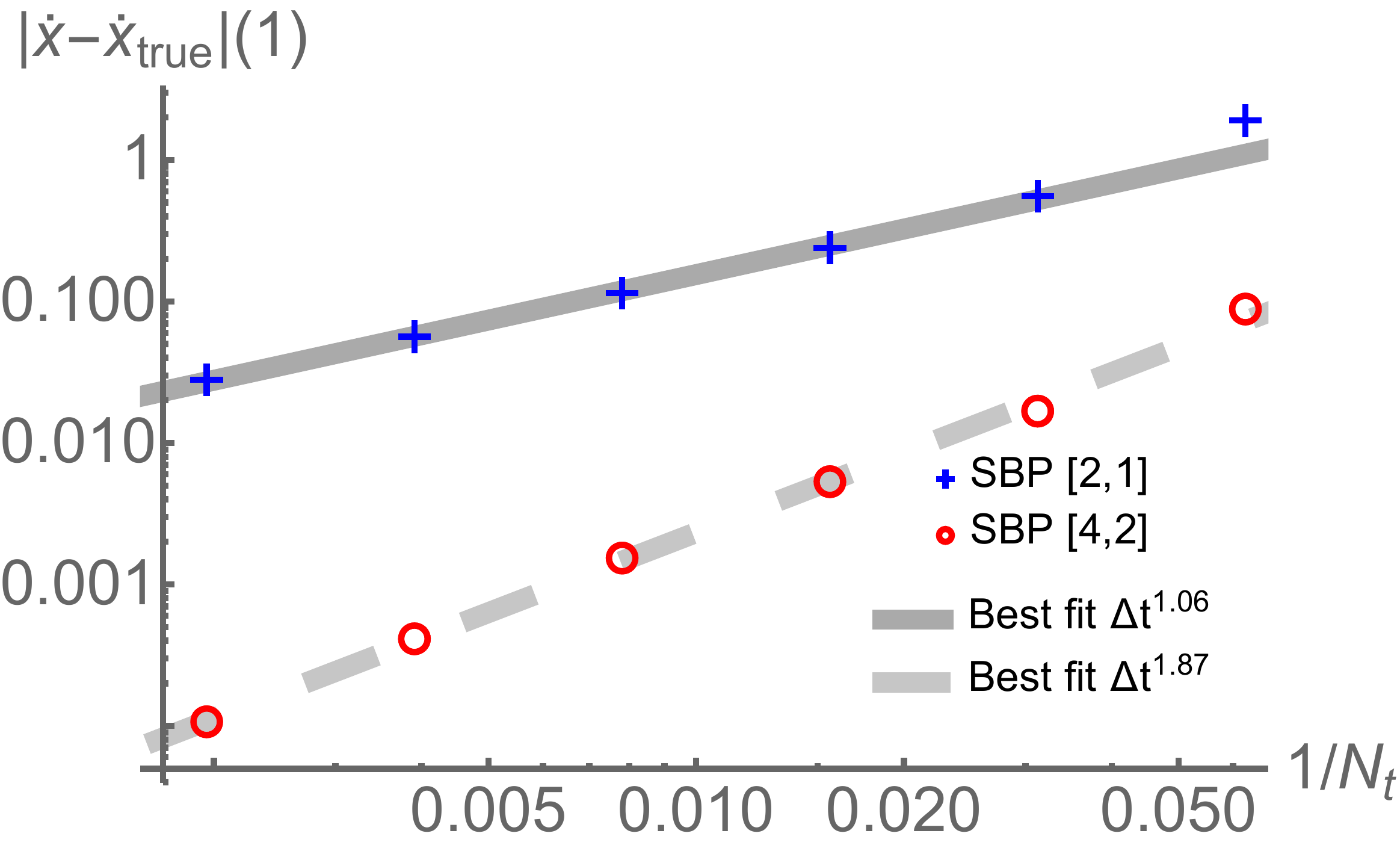}
    \caption{Deviation of the value (left plot) and the derivative (right plot) between the optimal numerical solution of \cref{eq:IVPnonlinfuncReg} and the true solution at the final time $t_2=1$. The deviation for the regularized \texttt{SBP21} operator is given as blue crosses with best fit $\Delta t^{2.12}$, the one for the regularized \texttt{SBP42} operator as red circles with $\Delta t^{3.35}$. The best power law fit to the data is shown as gray lines. Both orders of the solution improve steadily with diminishing $\Delta t\sim 1/N_t$. We find that our approach exhibits convergence in the values with one higher power in the grid spacing than for the derivative, where for \texttt{SBP21} we find as best fit $\Delta t^{1.06}$ and for \texttt{SBP42} we have $\Delta t^{1.87}$.}
    \label{fig:nonlfuncsol}
\end{figure}

Let us consider how the new discretization prescription performs on this second order non-linear problem. The continuum solution is given in terms of the Jacobi elliptic function and its inverse, i.e. it is not polynomial. Thus neither the \texttt{SBP21} nor the \texttt{SBP42} operator are able to reproduce it exactly. We find minute oscillations around the true solution exhibited by the data points in \cref{fig:classtrajnonlIVPSBPReg}, which diminish monotonously with grid refinement, as expected from a stable procedure. The deviation of the numerical solution ${\bf x}_1$ from the true solution is shown in the left panel of \cref{fig:nonlfuncsol} as blue crosses for the \texttt{SBP21} operator and as red circles for the \texttt{SBP42} operator.

One finds that compared to the linear system, the error made in the non-linear system is around one order of magnitude larger at the same lattice spacing for the \texttt{SBP21} operator. However the order of convergence remains close to quadratic with a best fit $\Delta t^{2.12}$ shown by the gray lines in \cref{fig:nonlfuncsol}. For the \texttt{SBP42} operator, we find that the error at $N_t=16$ improves by one order of magnitude and convergence to the continuum limit proceeds with $\Delta t ^{3.35}$. The expected behavior from solutions of ODEs \cite{svard2019convergence,svard_convergence_2021} in this case is $\Delta t ^{3}$, and while our method seems to perform better in this scenario, we believe that to be a coincidence.

In line with established results for the conventional numerical treatment of differential equations, we see in \cref{fig:nonlfuncsol} that the error in the derivative $\dot x=[\mathds{D}{\bf x}_1]({N_t})$ of the numerical solution shows convergence with one full order less than the values of the solution itself. As plotted in the right panel of \cref{fig:nonlfuncsol} we obtain for the regularized \texttt{SBP21} operator convergence for the derivative according to $\Delta t^{1.06}$, while the \texttt{SBP42} operator exhibits $\Delta t^{1.87}$.

Having established the applicability and convergence properties of our novel discretization approach for both a linear and non-linear second order differential equation of motion, let us continue to treat systems with equations of motion that feature different orders of time derivatives.

\section{Discretization of first order derivative terms}
\label{sec:generalizedVP}

So far we have considered the simplest case of physical systems with classical equations of motion that contain a second derivative in time. These follow naturally from the conventional formulation of the continuum variational principle, based on an action that is written in terms of a Lagrangian. As has been shown in \cite{galley_classical_2013}, the variational principle is able to accommodate a much larger variety of systems, including those with dissipative forces, which are not time-reversal invariant.
Such systems exhibit equations of motion, which contain also first order derivatives in time. The crucial step is to realize that, one may generalize the classical variational principle by adding to the Lagrangian $L$ another functional $\Lambda$ that may depend on both the forward and backward path\footnote{It is interesting to note that the term $\Lambda$ also arises naturally in the classical limit of the Schwinger-Keldysh contour formalism of the quantum path integral for dissipative systems, after integrating out the environment degrees of freedom. In that context it is known as the so called Feynman-Vernon influence functional.} and their derivatives as follows 
\begin{align}
 \nonumber &S_{\rm GIVP}[x_1(t),\dot x_1(t),x_2(t),\dot x_2(t)]=\\
    &\int_{t_1}^{t_2} dt\Big( {\cal L}[x_1(t),\dot x_1(t)] - {\cal L}[x_2(t),\dot x_2(t)] + \Lambda[x_1(t),\dot x_1(t),x_2(t),\dot x_2(t)]\Big).\label{eq:GIVPcont}
\end{align}
Ref. \cite{galley_classical_2013}, aided by ref.\cite{berges_quantum_2007}, shows in detail that the classical equations of motion also for this generalized variational principle are obtained by going over to relative $x_-=x_1-x_2$ and centered coordinates $x_+=(x_1+x_2)/2$ with the defining equation
\begin{align}
\left. \frac{\delta S_{\rm GIVP}[x_{\pm}]}{\delta x_-}\right|_{x_-=0,x_+=x_{\rm class}}=0. \label{eq:GIVPder}
\end{align}
The stability properties of these systems, as shown explicitly in ref.~\cite{sieberer2016keldysh}, can also be formulated in terms of a generalized Noether's theorem in which e.g. the time-dependence of the total system energy is correctly captured. This immediately invites us to apply the discretization prescription developed in the previous section to two systems often considered in the literature, the one which features the defining equation of the exponential function as equation of motion, as well as the damped harmonic oscillator. By considering these two examples, we acquire intuition in how to construct the appropriate continuum functional $\Lambda$ of \cref{eq:GIVPcont}, in order to describe systems, which feature a differential equation of motion also containing first order derivatives.

\subsection{A purely first order system}
\label{sec:firstordersys}
\begin{figure}
    \centering
    \includegraphics[scale=0.28]{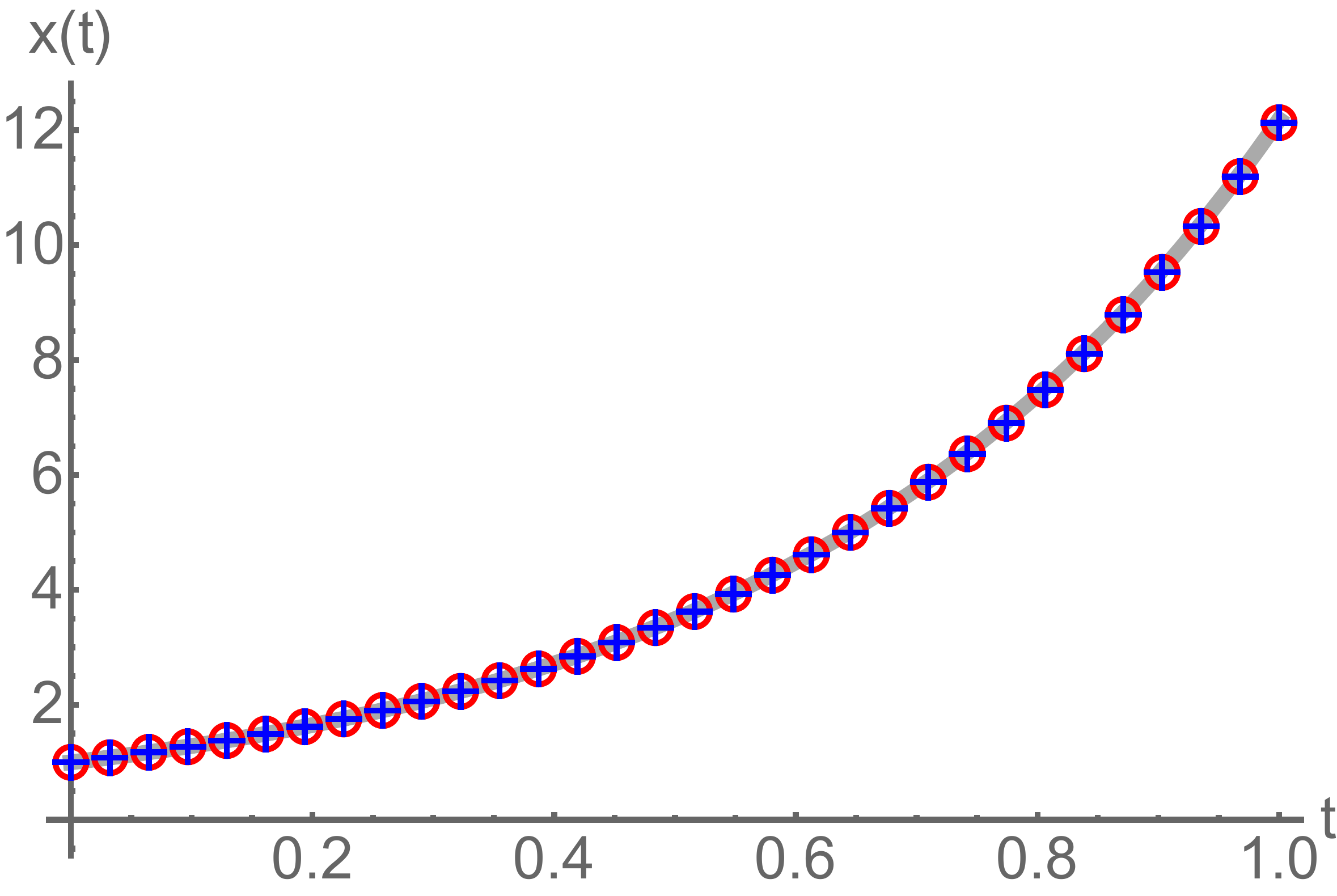}
    \caption{Numerical solution of the discretized path ${\bf x}_1$ (red circles) and ${\bf x}_1$ (blue crosses) ($N_t=32$) that optimize the functional \cref{eq:expfuncIVPfuncReg} corresponding to the discretized $S_{\rm GIVP}$. We use the regularized \texttt{SBP21} operator in time. Continuum solution of the Euler-Lagrange equation $\dot x(t)-\frac{5}{2} x(t)=0$ is shown as solid gray line.}
    \label{fig:expfuncSBPRegWeak}
\end{figure}

Our goal here is to determine the classical trajectory of a system, which features as its equation of motion the defining equation of the exponential function
\begin{align}
\dot x(t)-\kappa x(t)=0. \label{eq:expfuncdef}
\end{align}
Since there are no second order derivatives present in \cref{eq:expfuncdef} we do not need to supply the standard kinetic term to {\cal L} in \eqref{eq:GIVPcont}. The term linear in $x(t)$ can be thought of as arising from a potential contribution in ${\cal L}$, which must contain one power of $x_+$ and one power of $x_-$ similar to our argument in \cref{eq:constrsecord}. The new ingredient is the term that features a single time derivative. It has to emerge from $S_{\rm GIVP}$ after functional differentiation with respect to $x_-$. This behavior is achieved by choosing the following Lagrangian and $\Lambda$ functional
\begin{align}
{\cal L}=- \frac{1}{2}\kappa  x^2(t), \quad \Lambda=\dot x_+ x_-, \label{eq:expfunceom}
\end{align}
which amounts to the joint Lagrangian
\begin{align}
{\rm L}=  -\kappa x_+(t)x_-(t) + \dot x_+ x_-. \label{eq:expfunceomL}
\end{align}
Note that if one rewrites \cref{eq:expfunceomL} explicitly in terms of $x_1$ and $x_2$, the contribution from $\Lambda$ indeed does not factorize into terms that  depend on $x_1$ or $x_2$ separately. 

\begin{figure}
    \centering
    \includegraphics[scale=0.24]{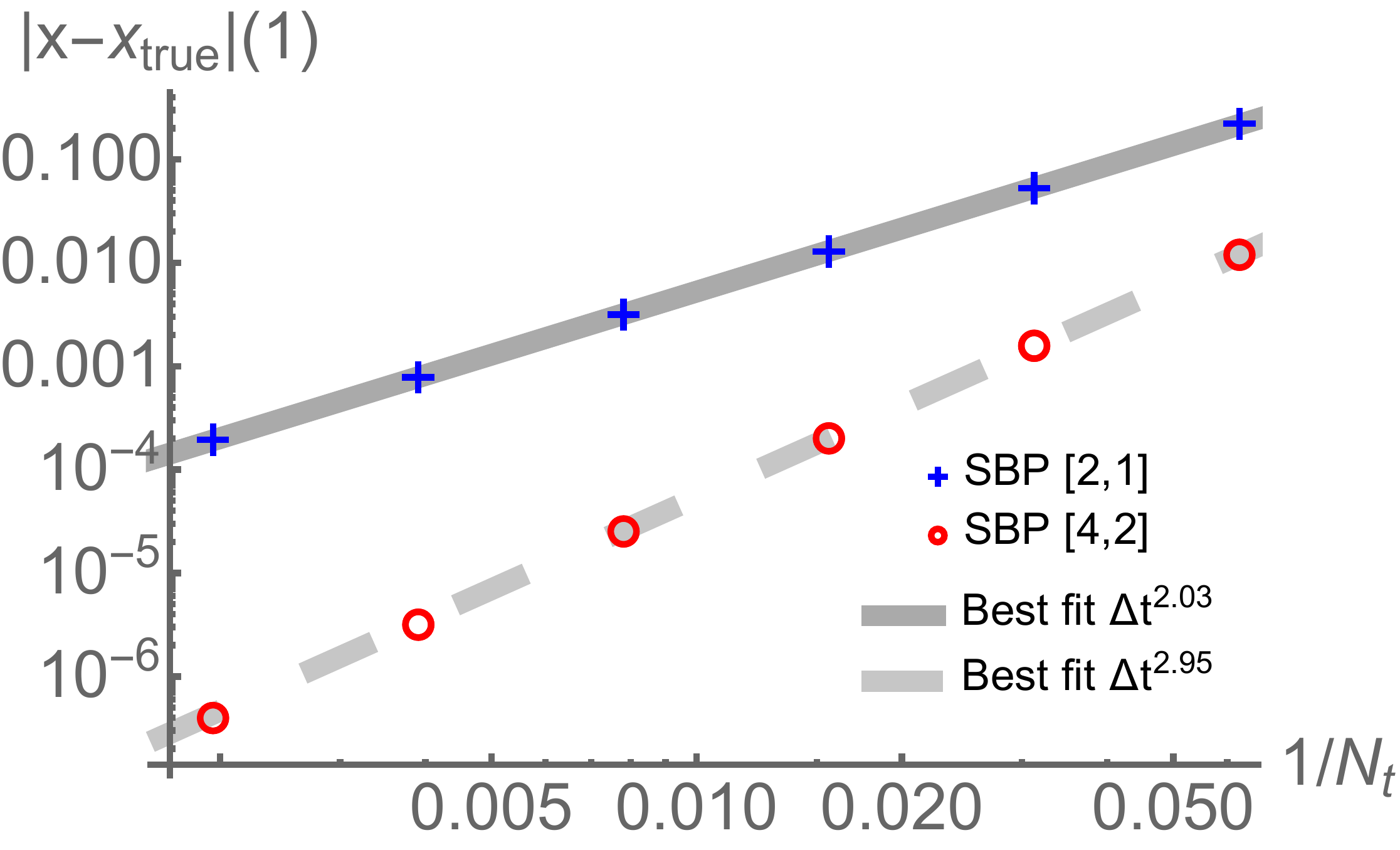}
    \includegraphics[scale=0.24]{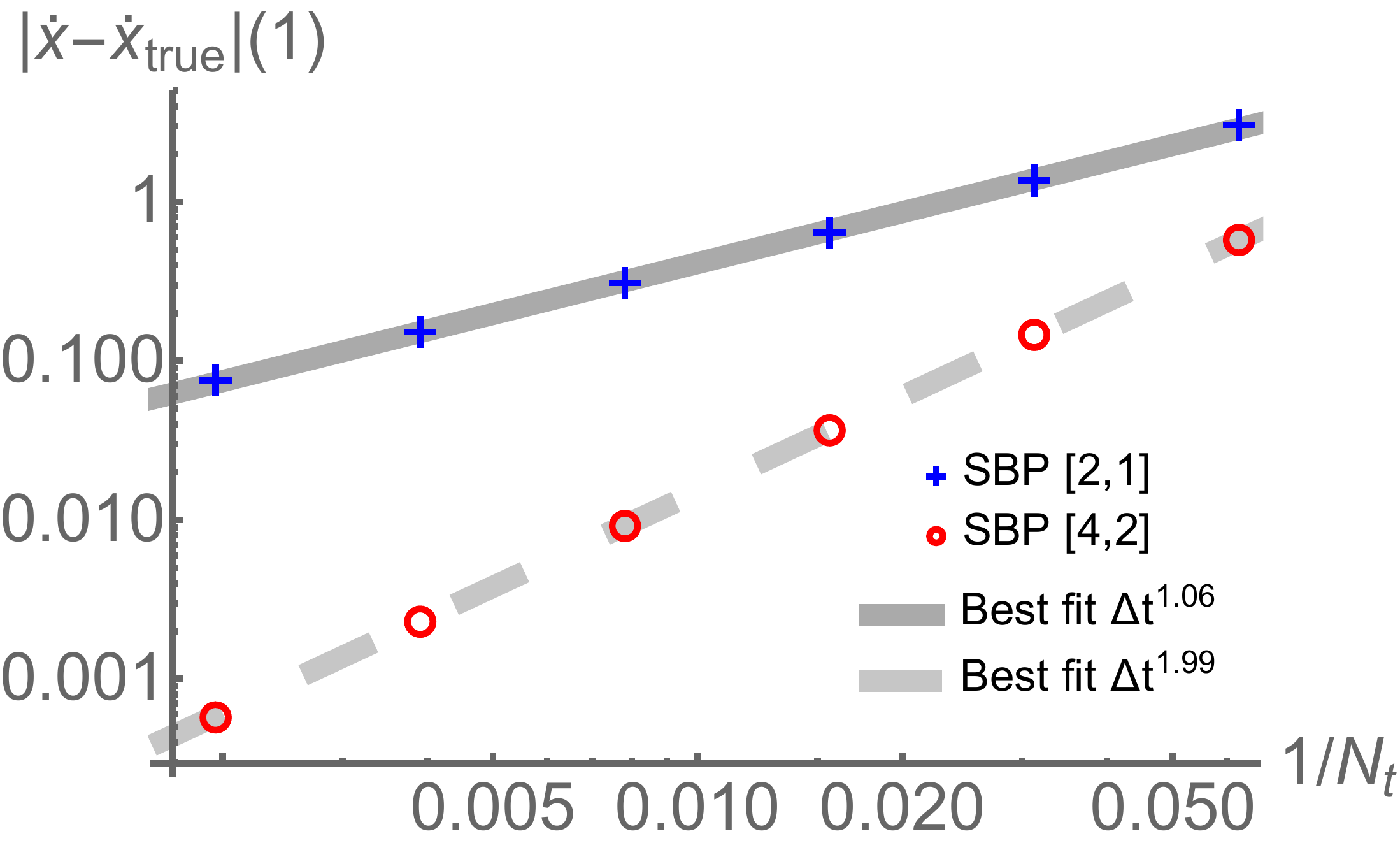}
    \caption{Deviation of the value (left plot) and the derivative (right plot) between the optimal numerical solution of \cref{eq:expfuncIVPfuncReg} and the true solution at the final time $t_2=1$. The deviation for the regularized \texttt{SBP21} operator is given as blue crosses with best fit $\Delta t^{2.03}$, the one for the regularized \texttt{SBP42} operator as red circles and $\Delta t^{2.95}$. The best power law fit to the data is shown as gray lines. Both orders of the solution improve steadily with diminishing $\Delta t\sim 1/N_t$. We find that our approach exhibits convergence in the values with one higher power in the grid spacing than for the derivative, where for \texttt{SBP21} we find $\Delta t^{1.06}$ and for \texttt{SBP42} we have $\Delta t^{1.99}$.}
    \label{fig:expfuncsol}
\end{figure}
Using the strategy developed in the previous section, let us write down the discretized action functional, keeping in mind that for a first order equation only the initial position needs to be supplied at the beginning of ${\bf x}_1$. Correspondingly only the position information needs to be matched at the end of the paths 
\begin{align}
\nonumber \mathds{S}_{\rm GIVP}&= \Big\{  -  \frac{1}{2}\kappa  ({\bf x}_1)^{\rm T} \mathds{H} {\bf x}_1 \Big\} -\Big\{  -  \frac{1}{2}\kappa ({\bf x}_2)^{\rm T} \mathds{H} {\bf x}_2 \Big\}\\
\nonumber &+ \frac{1}{2} \Big(\bar{\mathds{D}}\big({\bar{\bf x}}_1+{\bar{\bf x}}_2\big)\Big)^{\rm T} \bar{\mathds{H}} ({\bf x}_1-{\bf x}_2)  \\
 &+ \lambda_1 (x_1(0)-x_i) + \lambda_3 (x_1(N_t)-x_2(N_t)).\label{eq:expfuncIVPfuncReg}
\end{align}
The solutions ${\bf x}_1$ (red circles) and ${\bf x}_2$ (blue crosses) to this equation on the interval $t\in[0,1]$ discretized with $N_t=32$ equidistant steps $\Delta t$ and a $\kappa=\frac{5}{2}$ with initial condition $x_i=1$ produces the data shown in \cref{fig:expfuncSBPRegWeak}.

In \cref{fig:expfuncsol} we plot the difference between the classical trajectory obtained from \cref{eq:expfuncIVPfuncReg} and the true solution at the final time $t_2=1$ focussing on the value itself in the left plot and the derivative in the right plot. The results for the  regularized \texttt{SBP21} operator are given as blue crosses, those for the \texttt{SBP42} operator as red circles. Since the solution is not a simple polynomial, the \texttt{SBP42} operator cannot exactly integrate it. Both the \texttt{SBP21} and \texttt{SBP42} cases show the same convergence rates, as observed in the conventional formulation of IVPs (c.f. ref.~\cite{svard_convergence_2021}). We find again that the convergence is one order higher in the values of the solution than in the derivative of the solution. The \texttt{SBP21} operator yields a $\Delta t^{2.03}$ improvement for the values of $x_1(1)$, while the \texttt{SBP42} operator exhibits $\Delta t^{2.95}$.

\subsection{The damped harmonic oscillator}
\begin{figure}
    \centering
    \includegraphics[scale=0.28]{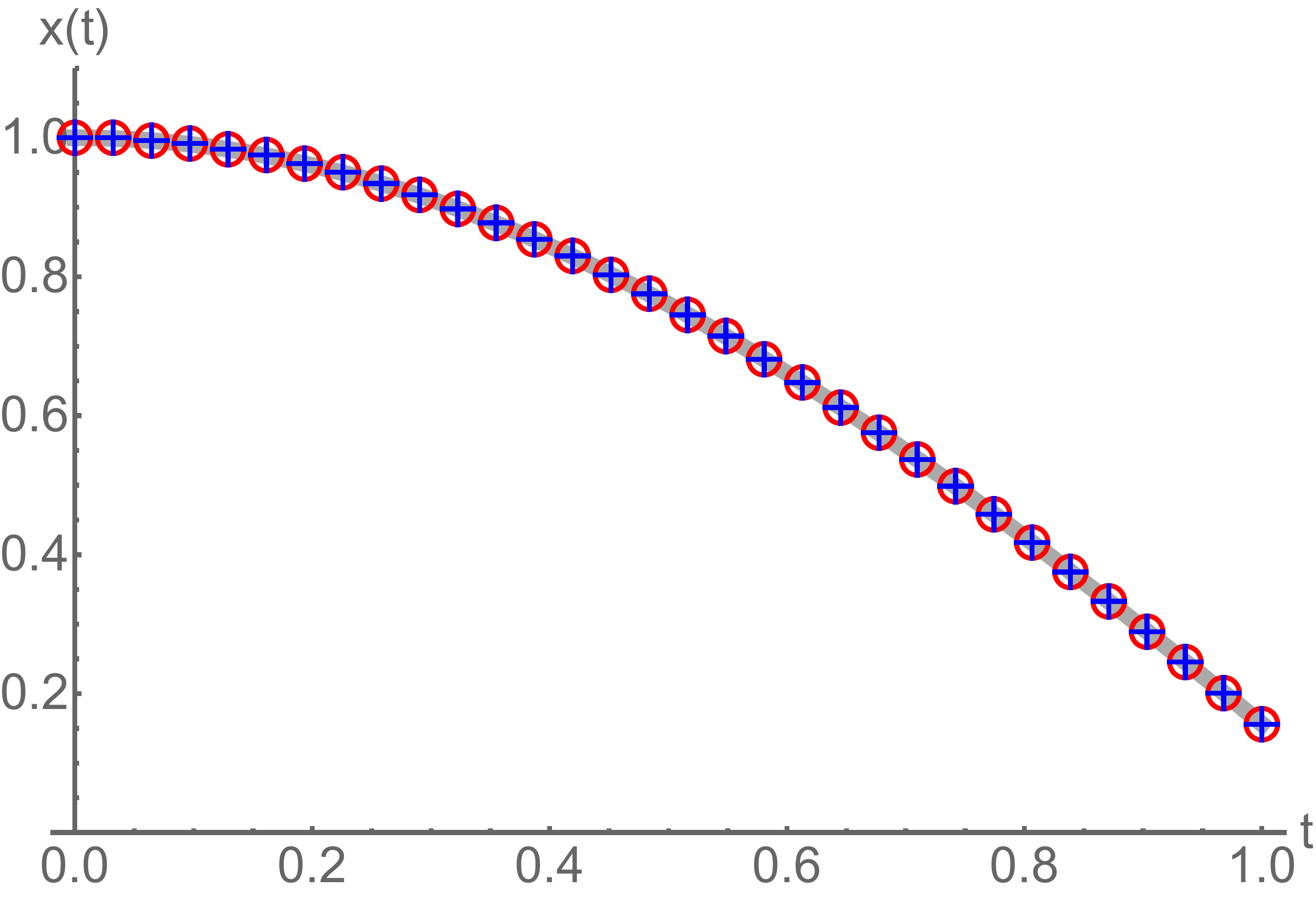}
    \caption{Numerical solution (red dots) of the discretized path ${\bf x}_1$ (red circles) and ${\bf x}_2$ (blue crosses) ($N_t=32$) that optimize the functional \cref{eq:damphoIVPfuncReg} corresponding to the discretized $S_{\rm GIVP}$. We use the regularized \texttt{SBP21} operator in time. Continuum solution of the Euler-Lagrange equation $\mu \ddot x(t)+\xi \dot x(t) + \kappa x(t)=0$ is shown as solid gray line.}
    \label{fig:dampHOSBPRegWeakShort}
\end{figure}

\begin{figure}
    \centering
    \includegraphics[scale=0.24]{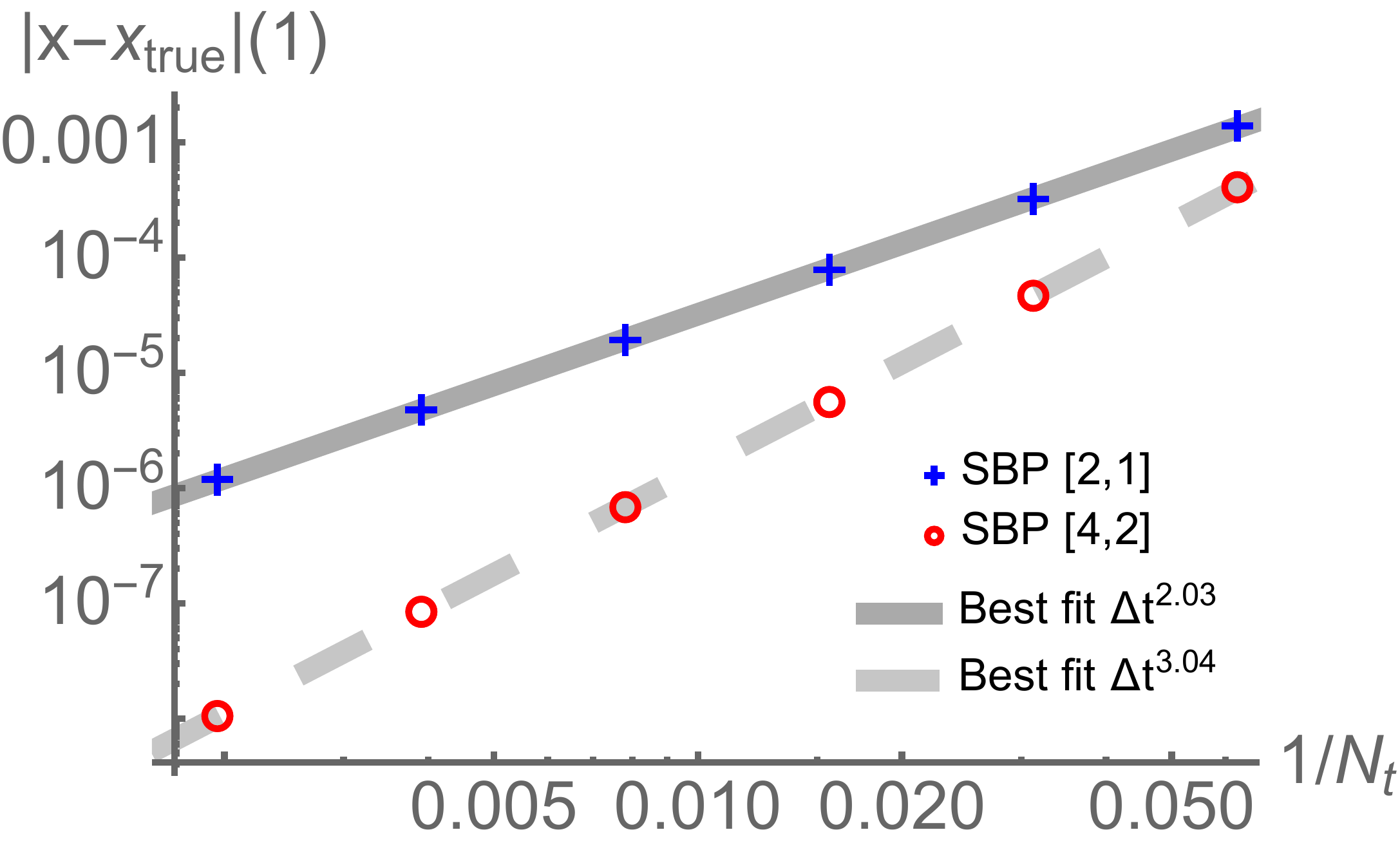}
    \includegraphics[scale=0.24]{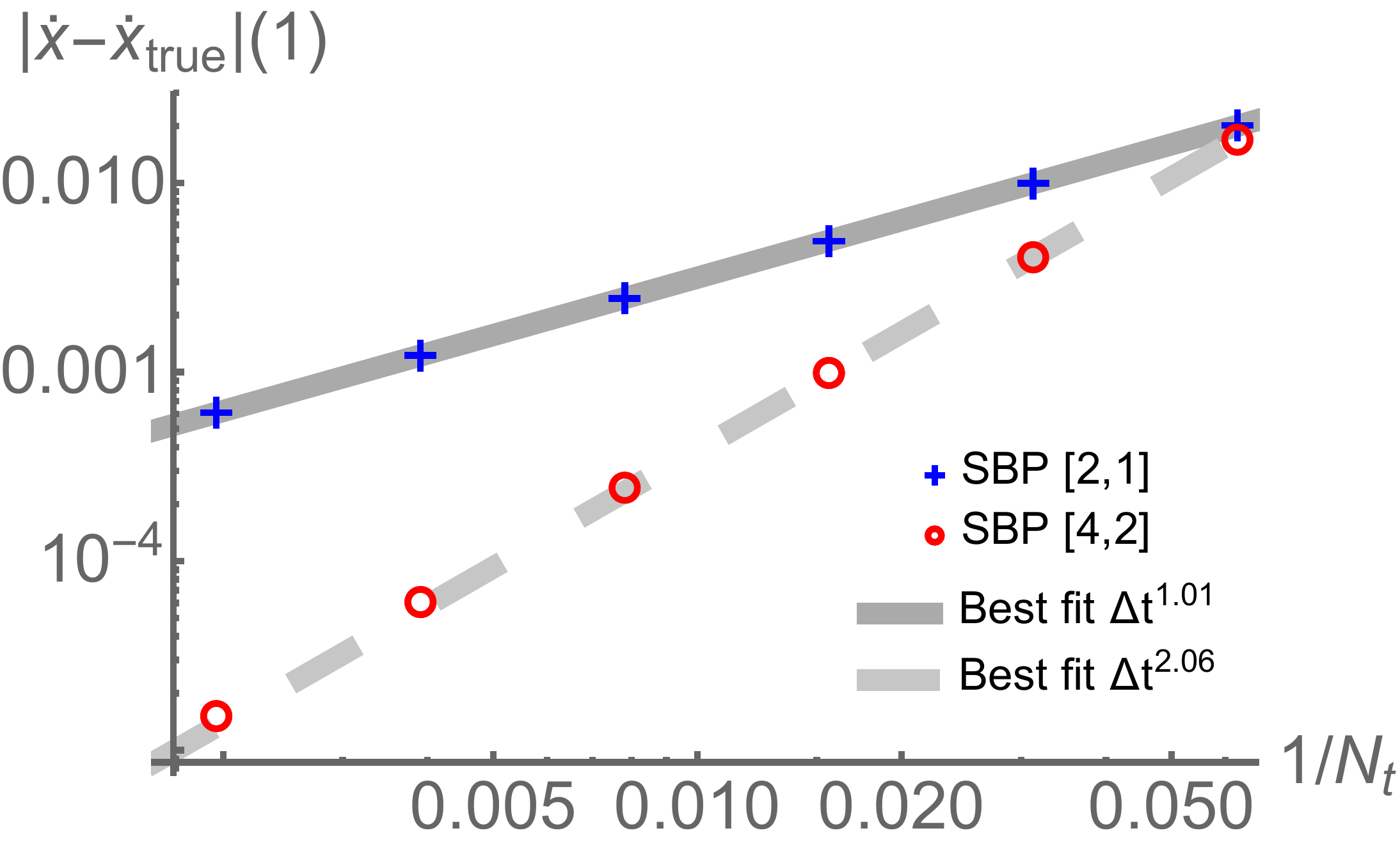}
    \caption{Deviation of the value (left plot) and the derivative (right plot) between the optimal numerical solution of \cref{eq:damphoIVPfuncReg} and the true solution at the final time $t_2=1$. The deviation for the regularized \texttt{SBP21} operator is given as blue crosses with best fit $\Delta t^{2.03}$, the one for the regularized \texttt{SBP42} operator as red circles with $\Delta t^{3.04}$. The best power law fit to the data is shown as gray lines. Both orders of the solution improve steadily with diminishing $\Delta t\sim 1/N_t$.  We find again that our approach exhibits convergence in the values with one higher power in the grid spacing than for the derivative, where for \texttt{SBP21} we find $\Delta t^{1.01}$ and for \texttt{SBP42} we have $\Delta t^{2.06}$.}
    \label{fig:damphofuncsol}
\end{figure}
As final challenge let us now turn to a physics system, which exhibits both first and second order derivatives in its equation of motion: the damped harmonic oscillator. The damped harmonic oscillator is characterized by an Euler-Lagrange equation that reads
\begin{align}
\mu \ddot x(t)+\xi \dot x(t) + \kappa x(t)=0. \label{eq:dampho}
\end{align}
Even though this system underlies a wealth of experimentally relevant phenomena, the conventional formulation of classical mechanics is unable to accommodate it in terms of a classical Lagrangian. In the generalized approach of ref.~\cite{galley_classical_2013} the kinetic and conservative force term are captured by the Lagrangian and the dissipative term is included via the $\Lambda$ functional. Remembering that only terms linear in $x_-$ survive the stationarity condition and using the intuition we built in the preceding sections, we can now write down the corresponding expressions for the functionals of the damped harmonic oscillator
\begin{align}
{\cal L}=\frac{1}{2} \mu \dot x^2(t) - \kappa x^2(t), \quad \Lambda=-\xi \dot x_+ x_- ,\label{eq:dmphocontact}
\end{align}
which correspond to the joint Lagrangian
\begin{align}
{\rm L}= \mu \dot x_+(t) \dot x_- - 2\kappa x_+(t)x_-(t) -\xi \dot x_+ x_- .\label{eq:dmphocontactL}
\end{align}
Inserting the above into \cref{eq:GIVPder} immediately yields \cref{eq:dampho}. The discretized joint action functional based on \cref{eq:dmphocontactL} reads 
\begin{align}
\nonumber \mathds{S}_{\rm GIVP}&= \Big\{  \frac{1}{2} \mu(\bar{\mathds{D}}{\bar{\bf x}}_1)^{\rm T} \bar{\mathds{H}} (\bar{\mathds{D}}\bar{{\bf x}}_1) -  \frac{1}{2} \kappa ({\bf x}_1)^{\rm T} \mathds{H} {\bf x}_1 \Big\} -\Big\{  \frac{1}{2} \mu(\bar{\mathds{D}}{\bar{\bf x}}_2)^{\rm T} \bar{\mathds{H}}(\bar{\mathds{D}}\bar{{\bf x}}_2) -  \frac{1}{2} \kappa ({\bf x}_2)^{\rm T} \mathds{H} {\bf x}_2 \Big\}\\
\nonumber &-\xi \frac{1}{2} \Big(\bar{\mathds{D}}\big({\bar{\bf x}}_1+{\bar{\bf x}}_2\big)\Big)^{\rm T} \bar{\mathds{H}} ({\bf x}_1-{\bf x}_2)  \\
\nonumber &+ \lambda_1 (x_1(0)-x_i) + \lambda_2((\mathds{D}{\bf x}_1)(0)-\dot x_i) \\
&+ \lambda_3 (x_1(N_t)-x_2(N_t)) + \lambda_4 ( (\mathds{D}{\bf x}_1)(N_t)-(\mathds{D}{\bf x}_2)(N_t) ).\label{eq:damphoIVPfuncReg}
\end{align}
Searching for the extremum of this functional numerically using the parameters $\mu=0.5$, $\kappa=1$, $\xi=0.00071$ (c.f. ref.~\cite{tsang_slimplectic_2015}) with initial conditions $x_i=1$ $\dot x_i=0$ on discretized paths with $N_t=32$ steps and regularized \texttt{SBP21} operator leads to the results for ${\bf x}_1$ (red circles) and ${\bf x}_2$ (blue crosses) shown in \cref{fig:dampHOSBPRegWeakShort}. 

\begin{figure}
    \includegraphics[scale=0.24]{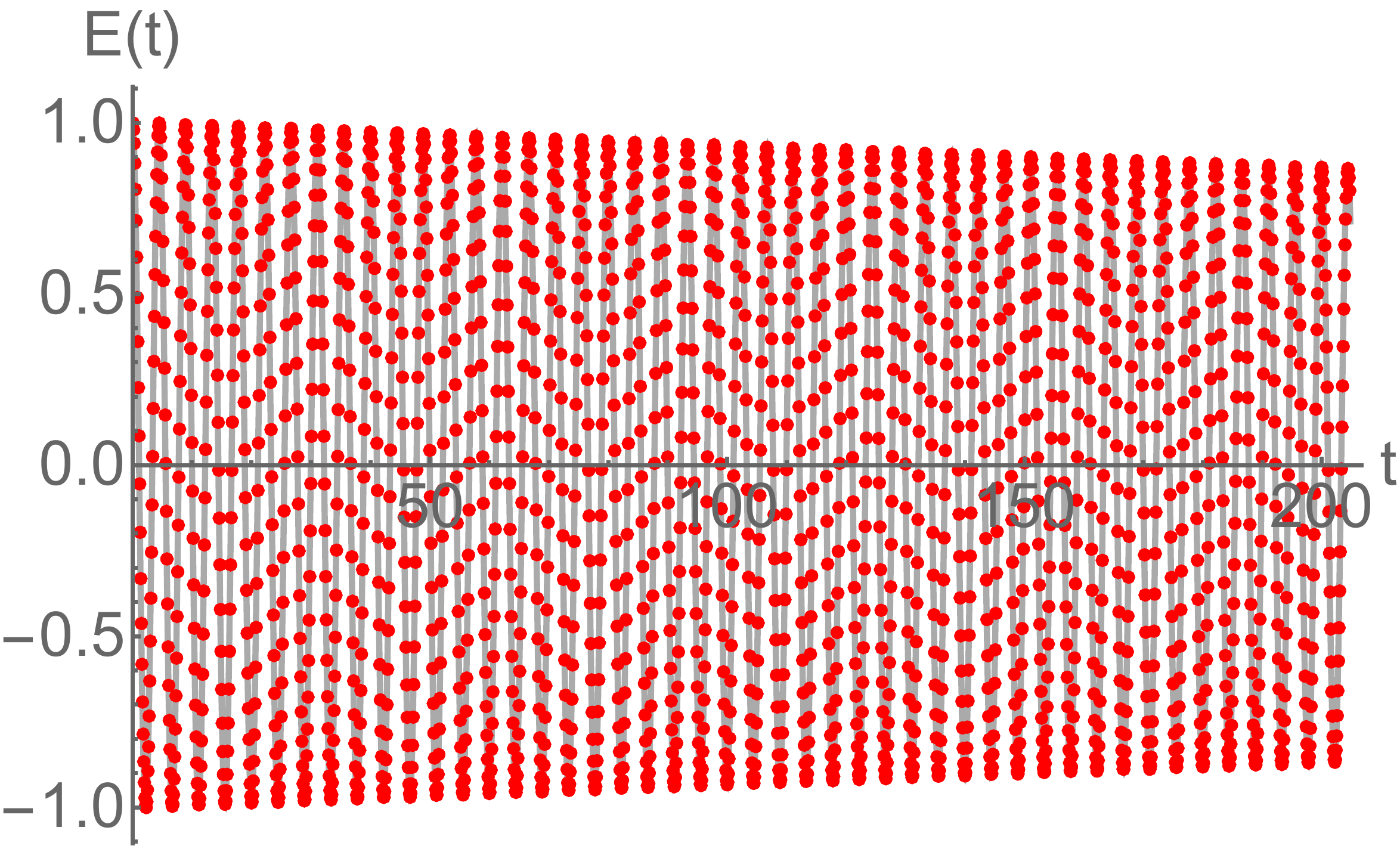}
    \includegraphics[scale=0.24]{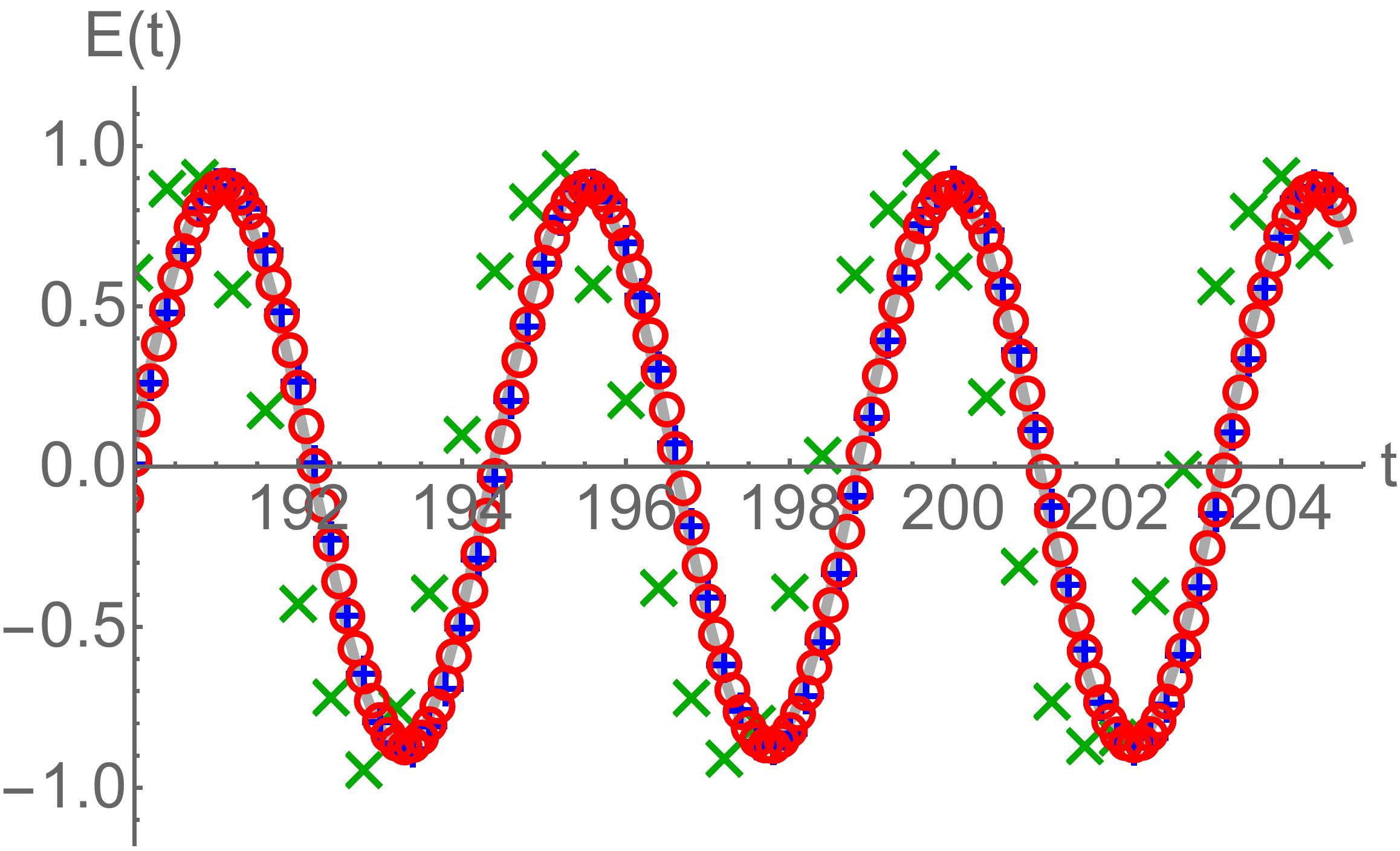}
    \caption{(Left panel) Numerical solution (red dots) of the discretized path ${\bf x}_1$ ($N_t=2048$, $\Delta t=0.1$) that optimizes the functional \cref{eq:damphoIVPfuncReg} corresponding to the discretized $S_{\rm GIVP}$ using the \texttt{SBP42} operator in time. The path  ${\bf x}_2$ takes on the same values. Continuum solution of the Euler-Lagrange equation $\mu \ddot x(t)+\xi \dot x(t) + \kappa x(t)=0$ is given as solid gray line. (Right panel) Comparison of the optimal solutions at a different number of grids points $N_t=512$ (green crosses), $N_t=1024$ (blue crosses) and $N_t=2048$ (red circles) at late times.}
    \label{fig:dampHOLatevals}
\end{figure}
In the left panel of \cref{fig:damphofuncsol} we present the deviation between the numerically determined critical point of \cref{eq:expfuncIVPfuncReg} and the true solution at the final time $t_2=1$. The deviation of the derivative is given in the right plot. Blue crosses denote the \texttt{SBP21} operator case, while red circles refer to the \texttt{SBP42} operator. The power law fits show that also for the damped harmonic oscillator the convergence order agrees with the expectations from the conventional formulation of IVPs \cite{svard_convergence_2021} and is one order higher in the values of the solution than in the derivative of the solution. The \texttt{SBP21} operator yields a $\Delta t^{2.03}$ improvement for the values of $x_1(1)$, while the \texttt{SBP42} operator exhibits $\Delta t^{3.04}$.

So far we have only investigated the short time behavior. However it is late-time stability that plays the most important role for the utility of a discretization scheme to the description of physical processes in practice. For initial boundary value problems this calls for so-called error-bounded schemes \cite{nordstrom_long_2018,kopriva_error_2017,nordstrom_error_2008}. In the literature the damped harmonic oscillator is often used as a non-trivial test-bed to evaluate the late-time stability and accuracy of numerical solvers. Let us therefore determine the numerical solution to \cref{eq:damphoIVPfuncReg} up to $t_2=204.8$, which is shown in the left panel of \cref{fig:dampHOLatevals}, based on the regularized \texttt{SBP42} operator at $\Delta t=0.1$. This choice of $t_2$ allows the system to pass through multiple oscillations and to show a visible reduction of the oscillation amplitude. As our SBP in time discretization is inherently implicit, we find numerically that while the solution degrades in accuracy as we increase the grid spacing $\Delta t\in\{0.1,0.133,0.2,0.4\}$ it remains bounded for all times. The behavior of the discrete solution for different grid spacings $\Delta t$ at late times, is shown in the right panel of \cref{fig:dampHOLatevals}. One finds that the solution converges as the grid spacing is decreased. The most pertinent error introduced by the discretization procedure appears to be an artificial phase shift, the dispersion error, which however vanishes as the continuum limit $\Delta t\to0$ is approached.

\begin{figure}
    \centering
    \includegraphics[scale=0.28]{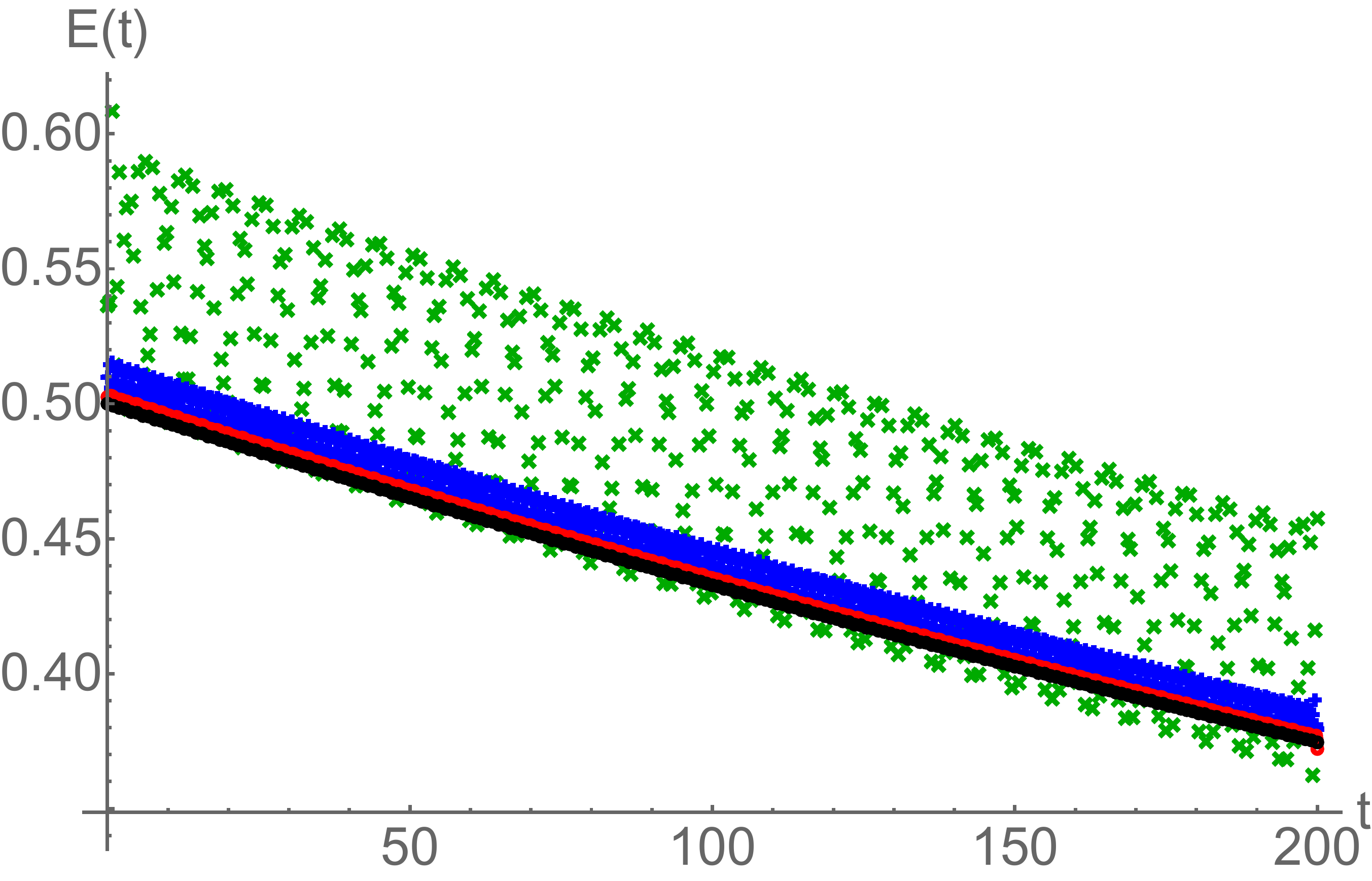}
    \caption{Numerical solution of the energy of the solution path ${\bf x}_1$ from functional \cref{eq:damphoIVPfuncReg} using the regularized \texttt{SBP42} operator. The energy is plotted at a different number of grids points $N_t=512$ (green crosses), $N_t=1024$ (blue crosses) and $N_t=2048$ (red circles) and compared to the continuum solution from the Euler-Lagrange equation, which is shown as solid black line. Note that refining the grid consistently leads to smaller oscillations and the numerical solution approaches the continuum result.}
    \label{fig:dampHOSBPReNrg}
\end{figure}

A common quality criterion for numerical solvers in the physical sciences is the reproduction of the system energy. We consider here as energy the following Hamiltonian ${\cal H}=T+V=\frac{1}{2}\mu\dot x^2(t)+\kappa x^2(t)$, which is plotted for different grid spacings in \cref{fig:dampHOSBPReNrg}. One can clearly see that the discretization procedure leads to an overall shift in the value of the energy and the appearance of oscillations around the mean value, known as dispersion and diffusion errors (for a detailed exposition of the dispersion errors of SBP operators see e.g.  ref.~\cite{linders2015uniformly} and references therein). However, both the shift, as well as the oscillations vanish with grid refinement and no artificial energy deviation or instability is observed. 

\section{Summary}
\label{sec:summary}

We have presented a new and unified discretization strategy for a wide range of initial value problems, formulated in terms of a generalized variational problem. The classical trajectory is obtained as the critical point of an action functional with doubled degrees of freedom, without the need to derive equations of motion. Our approach is based on summation-by-parts finite difference operators, which are regularized through a penalty term, associated with the initial conditions using affine coordinates. We introduce the initial conditions, as well as the identification of the forward and backward path in the IVP in \cref{sec:strongIVP} using Lagrange multipliers. The approach has been successfully applied to systems with purely second order, first order and mixed time derivatives in their equations of motion. Explicit scaling tests for each of these systems confirmed that our prescription converges to the true solution under grid refinement, with the same convergence rates, as expected from the conventional formulation of IVPs. Convergence in the values of the solution of the systems investigated here occurs with one higher order in the grid spacing $\Delta t$ compared to the derivative of the solution. Our approach is based on a single realization of the first order SBP finite difference operator and does not require a separate treatment of systems with equations of motion containing first or second order derivatives in time. 

Our study describes a genuinely novel and versatile approach to the computation of trajectories of classical systems without the need to resort to their equation of motion. The extension to higher dimensions for the treatment of partial differential equations is the natural next step and is work in progress. We are looking forward to applying the formalism to the discretization of theories with intrinsic constraints, in particular gauge theories, such as classical electromagnetism and Yang-Mills theory.

\section*{Acknowledgements}
A.~R. thanks Will Horowitz for stimulating discussions and gladly acknowledges support by the Research Council of Norway under the FRIPRO Young Research Talent grant 286883. J.~N. was supported by the Swedish Research Council grant nr. 2018-05084 and 2021-05484. The study has benefitted from computing resources provided by  
UNINETT Sigma2 - the National Infrastructure for High Performance Computing and Data Storage in Norway under project NN9578K-QCDrtX "Real-time dynamics of nuclear matter under extreme conditions"


\FloatBarrier

\begin{backmatter}

\section*{Competing interests}
  The authors declare that they have no competing interests.

\section*{Author's contributions}
    \begin{itemize}
         \item A. Rothkopf: project inception, development of the discretization and regularization prescription, scaling tests, writing and editing
         \item J. Nordstr\"om: guidance on the SBP-SAT formalism, development of the discretization formalism, literature review, writing and editing
    \end{itemize}

\appendix

\section{Noether Theorem for generalized variational problems}
\label{sec:Noether}

Following ref.~\cite{sieberer2016keldysh}  we retrace here how Noether's theorem connects global symmetries of the system with conserved quantities, the so called Noether charges in the case of doubled degrees of freedom. Note that transformations that act the same on the forward and backward path, due to the construction of $L$, will lead to vanishing Noether currents. To identify a finite conserved quantity, we need to consider transformations that act differently on the paths $x_{1,2}$. The relevant symmetry in our case is the invariance under opposite time translations. To formalize this statement, one converts an infinitesimal time shift $\tau$ in the argument of the paths via a Taylor expansion to ${\cal T}_\tau x_{1,2}(t)=x_{1,2}(t) \pm \tau \dot x_{1,2}(t) + {\cal O}(\tau^2)$. The Lagrangian ${\cal L}$ is a scalar, just as the paths and therefore transforms the same way as $\delta {\cal L}[x_{1,2}]/\delta \tau = \pm d{\cal L}[x_{1,2}]/dt $. One thus obtains
\begin{align}
\delta L&= \Big(\frac{\partial {\cal L}}{\partial x_1}\delta x_1 + \frac{\partial {\cal L}}{\partial \dot x_1}\delta \dot x_1 - \frac{\partial {\cal L}}{\partial x_2}\delta x_2 - \frac{\partial {\cal L}}{\partial \dot x_2}\delta \dot x_2 \Big)\\
\nonumber&= \Big(\frac{\partial {\cal L}}{\partial x_1}\delta x_1 - \frac{d}{dt}\frac{\partial {\cal L}}{\partial \dot x_1}\delta x_1 - \frac{\partial {\cal L}}{\partial x_2}\delta x_2 + \frac{d}{dt}\frac{\partial {\cal L}}{\partial \dot x_2}\delta x_2 + \frac{d}{dt} \Big\{ \frac{\partial {\cal L}}{\partial \dot x_1}\delta x_1 - \frac{\partial {\cal L}}{\partial \dot x_2}\delta x_2 \Big\}\Big)
\end{align}
where we have used integration by parts to arrive at the second line. The first four terms are nothing but the Euler Lagrange equations we obtained in \cref{eq:equivEL}, which vanish identically, so that we are left with
\begin{align}
\frac{\delta L}{\delta \tau}= \frac{d}{dt} \Big( \frac{\partial {\cal L}}{\partial \dot x_1}\dot x_1 + \frac{\partial {\cal L}}{\partial \dot x_2}\dot x_2\Big) = \frac{d}{dt}\Big( {\cal L}[x_1,\dot x_1]  + {\cal L}[x_2,\dot x_2] \Big)
\end{align}
The last equal sign arises from the fact that the Lagrangian itself transforms as a scalar. We thus arrive at the final expression for the conserved quantity as
\begin{align}
\frac{d}{dt}\Big( \pi_1 \dot x_1 + \pi_2\dot x_2 - {\cal L}[x_1,\dot x_1] - {\cal L}[x_2,\dot x_2] \Big)=\frac{d}{dt}\Big( H_+ + H_- \Big)=0\label{eq:NoetherEnergy}
\end{align}
This equation states that the total Hamiltonian, i.e. the total energy of the forward and backward path degrees of freedom is conserved in time. Had we instead started with time translations that acted the same on the forward and backward contour ${\cal T^\prime}_\tau x_{1,2}(t)=x_{1,2}(t) + \tau \dot x_{1,2}(t) + {\cal O}(\tau^2)$, the corresponding Noether charge would be the difference between the energies on the forward and backward path, i.e. 
\begin{align}
\frac{d}{dt}\Big( \pi_1\dot x_1 - \pi_2\dot x_2 - {\cal L}[x_1,\dot x_1] + {\cal L}[x_2,\dot x_2] \Big)=\frac{d}{dt}\Big( H_+ - H_- \Big)=0\label{eq:NoetherEnergyDiff}
\end{align}
telling us that not only is the total energy preserved but the energy difference between the forward and backward path must remain the same over time.

It is interesting to realize that the variational principle of \cref{eq:varIVPdefEq}, which has been derived here in a fully classical context, identically arises as the classical limit of IVPs in quantum field theory. As discussed in detail in ref.~\cite{berges_quantum_2007}, formulating initial value problems in the language of Feynman's path integral necessitates the introduction of doubled degrees of freedom. The combination of the forward and backward path are referred to as the Schwinger-Keldysh time contour. It turns out that the relative path $x_-$ is related to the quantum contributions and taking $x_-\to0$ is therefore intimately related to the classical limit. The classical limit of taking $\hbar\to0$ actually enforces $x_-=0$. Ref.~\cite{berges_quantum_2007} shows that the variation of the joint action with respect to $x_-$ is the relevant expression that describes how the fluctuating quantum paths collapse onto the deterministic classical trajectory, which indeed emerges after taking the limit $x_-\to0$. In the context of Noether's theorem, as discussed in ref.~\cite{sieberer2016keldysh}, transformations that act equally on forward and backward path are associated with quantum Noether currents, which do not have a finite expectation value. On the other hand transformations that couple the forward and backward path can be considered as quantum transformations, which lead to classical Noether currents that in turn can have a finite expectation value even in the classical limit, as we saw in \cref{eq:NoetherEnergy}.


\bibliographystyle{stavanger-mathphys}


\bibliography{references}


\end{backmatter}


\end{document}